\documentclass[12pt]{amsart}

\usepackage{a4wide}
\usepackage{amsmath}
\usepackage{epsfig}
\usepackage{url}
\usepackage{enumerate}

\newcommand{\numberthis}{\addtocounter{equation}{1}\tag{\theequation}}

\newcommand{\Sc}{\mathcal{S}}
\newcommand{\Z}{\mathbb{Z}}
\newcommand{\K}{\mathbb{K}}
\newcommand{\F}{\mathbb{F}}
\newcommand{\C}{\mathbb{C}}
\renewcommand{\H}{\mathbb{H}}
\newcommand{\Q}{\mathbb{Q}}

\newcommand{\Rb}{\tilde{R}}
\renewcommand{\P}{\mathbb{P}}
\newcommand{\CP}{\mathbb{P}}

\newcommand{\PU}{\rm PU}
\newcommand{\tr}{\rm tr}

\usepackage{chngcntr}
\counterwithin{table}{section}
\counterwithin{figure}{section}

\newtheorem{prop}{Proposition}[section]

\newtheorem{thmintro}{Theorem}
\newtheorem{thm}{Theorem}[section]
\def\cqfd{\mbox{}\nolinebreak\hfill$\Box$\medbreak\par}
\newenvironment{pf}{\noindent\textbf{Proof:}}{\cqfd}
\theoremstyle{definition}
\newtheorem{dfn}{Definition}[section]
\newtheorem{rk}{Remark}[section]

\newcommand{\quot}[2]{#2\backslash#1}
\newcommand{\quotr}[2]{#1/#2}

\newcommand{\gr}{\nabla}

\newcommand{\di}{\displaystyle}

\usepackage[matrix,arrow]{xy}
\xyoption{all}

\date{Jan 17, 2017}

\usepackage{array}

\title{Non-arithmetic lattices and the Klein quartic}

\author{Martin Deraux}

\address{Martin Deraux, Institut Fourier, Universit\'e Grenoble Alpes, France}
\email{martin.deraux@univ-grenoble-alpes.fr}

\begin{document}

\begin{abstract}
  We give an algebro-geometric construction of some of the
  non-arithmetic ball quotients constructed by the author, Parker and
  Paupert. The new construction reveals a relationship
  between the corresponding orbifold fundamental groups and the
  automorphism group of the Klein quartic, and also with groups
  constructed by Barthel-Hirzebruch-H\"ofer and
  Couwenberg-Heckman-Looijenga.
\end{abstract}

\maketitle

\section{Introduction}\label{sec:intro}

It is by now well known that there exist non-arithmetic lattices in
$\PU(2,1)$, the group of biholomorphisms of the complex hyperbolic
plane $\H^2$. The first examples were due to
Mostow~\cite{mostowpacific}, who gave explicit sets of matrices that
generate non-arithmetic lattices (the fact that these groups are
indeed lattices was proved by constructing explicit fundamental
domains for their action on the complex hyperbolic plane).

A similar approach (using different kinds of fundamental domains) was
used by the author in joint work with Parker and Paupert~\cite{dpp2},
and this allowed us to increase the number of known commensurability
classes of non-arithmetic lattices in $\PU(2,1)$ (there are currently
22 known classes, see~\cite{thealgo}).

The main goal of the present paper is to give an alternate
construction of the lattices $\Sc(p,\overline{\sigma}_4)$
of~\cite{dpp2}, that avoids fundamental domains altogether. For a
description of these groups, see section~\ref{sec:s4c}
(or~\cite{thealgo} for much more detail). Recall that, for
$p=4,5,6,8,12$, the groups $\Sc(p,\overline{\sigma}_4)$ give pairwise
incommensurable non-arithmetic lattices, that are not commensurable to
any Deligne-Mostow lattice.

Our results give a connection between these lattices and the
automorphism group of the Klein quartic, which is the closed Riemann
surface of genus 3 with largest possible automorphism group, i.e. such
that equality holds in the Hurwitz bound $|\textrm{Aut}|\leq
84(g-1)=168$.  It is well known that the corresponding automorphism
group is the unique simple group of order 168, and that it is also
isomorphic to the general linear group $GL_3(\F_2)$ in three variables
over the finite field with two elements. This group also turns out to
be the projectivization of a unitary group generated by complex
reflections, which appears as the group $G_{24}$ in the Shephard-Todd
list. The quotient $\quotr{\P^2}{G_{24}}$ is isomorphic as a normal
complex space to the weighted projective plane $X=\P(2,3,7)$ (this can
be checked by deriving invariant polynomials from the equation of the
Klein quartic, see section~\ref{sec:weighted} for details). The branch
locus of the corresponding quotient map $\P^2\rightarrow
\quotr{\P^2}{G_{24}}$ is an irreducible curve $M$ of degree 21 in
$\P(2,3,7)$, which we refer to as the Klein discriminant. This can be
rephrased to the statement that the pair $(X,(1-\frac{1}{2})M)$ is a
complex projective orbifold, i.e. an orbifold uniformized by $\P^2$,
which we refer to as the Klein orbifold.

We will see that the pair $(X,(1-\frac{1}{3})M)$ obtained from the
Klein orbifold by changing the multiplicity of $M$ from 2 to 3 is also
an orbifold, but now uniformized by the complex hyperbolic plane
$\H^2$. More generally, our main result will state that, for
$p=3,4,5,6,8,12,\infty$, the pair $(X,(1-\frac{1}{p})M)$ \emph{can be
  modified} to a complex hyperbolic orbifold. Indeed, for $p>3$,
changing the multiplicity of the ramification divisor is not enough,
much more drastic modifications are needed.

In order to explain these modifications, we start by establishing
notation. Let us denote by $\mathcal{H}$ the line arrangement in
$\P^2$ given by the mirrors of complex reflections in $G_{24}$, which
has 28 triple points and 21 quadruple points. The action of $G_{24}$
is transitive on the set of triple intersections, and also on the set
of quadruple intersections. We denote by $s_3$ (resp. $s_4$) the image
in $\P(2,3,7)$ of any triple (resp. quadruple) intersection.

Denote by $\widehat{Y}$ the blow-up of $\P^2$ at the 21 quadruple
intersections of $\mathcal{H}$. Similarly, $\widehat{Z}$ denotes the
blow-up at all 49 (triple or quadruple) singular points of
$\mathcal{H}$. We denote by $X=\quotr{\P^2}{G_{24}}=\P(2,3,7)$,
$Y=\quotr{\widehat{Y}}{G_{24}}$ and $Z=\quotr{\widehat{Z}}{G_{24}}$. With a slight abuse
of notation, we denote by $M$ either the Klein discriminant in $X$ or
its strict transform in $Y$ or $Z$. Similarly, $E$ denotes the
exceptional locus above $s_4$ either in $Y$ or in $Z$, and $F$ denotes
the exceptional locus in $Z$ above $s_3$.

We consider the pairs $(X^{(p)},D^{(p)})$ given in
Table~\ref{tab:pairs}.
\begin{thmintro}\label{thm:compact}
  For $p=3,5,8$ and $12$, the pairs $(X^{(p)},D^{(p)})$ are compact
  complex hyperbolic orbifolds. In other words, there exists a uniform
  lattice $G^{(p)}$ such that $\quot{\H^2}{G^{(p)}}$ is isomorphic
  to $X^{(p)}$ as a normal complex space, and the ramification divisor
  of the quotient map is given by $D^{(p)}$.
\end{thmintro}
For $p=4$ and $6$, the singularities of $(X^{(p)},D^{(p)})$ fail to be
log-terminal, so one cannot hope for these to be orbifold pairs (see
section~\ref{sec:logBMY}). In fact, the pair $(X^{(4)},D^{(4)})$ is
log-terminal at every point apart from $s_4$ (see
Proposition~\ref{prop:terminal_locus}), so we consider
$X_0^{(4)}=X\setminus\{s_4\}$. Similarly, $(X^{(6)},D^{(6)})$ is
log-terminal away from $s_3$, so we consider
$X_0^{(6)}=Y\setminus\{s_3\}$. One more pair can be used in order to
construct ball quotients, namely $(X^{(\infty)},D^{(\infty)})$. In
that case we take $X_0^{(\infty)}=Z\setminus M$. In all three cases
($p=4,6$ and $\infty$), we write $D_0^{(p)}=X_0^{(p)}\cap D$. We will
prove the following.
\begin{thmintro}\label{thm:noncompact}
  For $p=4$, $6$, and $\infty$ the pair $(X_0^{(p)},D_0^{(p)})$ is a
  non-compact complex hyperbolic orbifold of finite volume with one
  end. In other words, there exists a 1-cusped lattice $G^{(p)}$ such
  that $\quot{\H^2}{G^{(p)}}$ is isomorphic to $X_0^{(p)}$ as a
  normal complex space, and the ramification divisor of the quotient
  map is given by $D_0^{(p)}$.
\end{thmintro}
Theorems~\ref{thm:compact} and~\ref{thm:noncompact} will be proved by
showing that equality holds in the logarithmic Bogomolov-Miyaoka-Yau
inequality. More precisely, we use a result due to Kobayashi, Nakamura
and Sakai (see Theorem~\ref{thm:kns}), and a significant part of our
paper is devoted to checking that their result applies (the fact that
the hypotheses are indeed satisfied is our
Proposition~\ref{prop:verif}).

From the proof of the Kobayashi-Nakamura-Sakai result, it may seem
difficult to gather arithmetic information about the groups
$G^{(p)}$. Indeed, the existence is shown by solving the appropriate
Monge-Amp\`ere equation to produce a K\"ahler-Einstein metric with
negative Einstein constant. Under the assumption that equality holds
in the Bogomolov-Miyaoka-Yau inequality, one shows that the
K\"ahler-Eistein metric actually has constant holomorphic sectional
curvature. From there, it is not at all obvious how to obtain an
explicit description of the corresponding lattices in $U(2,1)$, say in
terms of explicit matrix generators.

In order to identify the lattices, we will use a description due to
Naruki of the fundamental group of the complement of the Klein
configuration of mirrors in $\P^2$, see~\cite{naruki}. This gives
enough explicit information about the orbifold fundamental group of
$(X^{(p)},D^{(p)})$ to show the following (see
section~\ref{sec:isomorphism}).
\begin{thmintro}\label{thm:iso}
  For $p=3,4,5,6,8$ or $12$, the lattices $G^{(p)}$ are conjugate in
  $\PU(2,1)$ to the lattices $\Sc(p,\overline{\sigma}_4)$. In
  particular, the lattices $G^{(p)}$, $p=4,5,6,8,12$,  are pairwise
  incommensurable non-arithmetic lattices, not commensurable to any
  Deligne-Mostow lattice.
\end{thmintro}
The case $p=\infty$ turns out to be less interesting, since the
corresponding lattice is arithmetic, see Proposition~\ref{prop:infty}.
The proof of Theorem~\ref{thm:iso} relies in part on previous joint
work of the author with Parker and Paupert, namely some of the
non-discreteness results in~\cite{dpp} (which we review in
section~\ref{sec:s4c}).

Putting together Theorems~\ref{thm:compact},~\ref{thm:noncompact}
and~\ref{thm:iso}, we get a new proof that the groups $\Sc(p,\bar
\sigma_4)$ are discrete, a fact which was known so far only by using
heavy computer work, see~\cite{dpp2} (and also~\cite{thealgo}).

The construction of ball quotients via uniformization (exploiting the
equality case in the Bogomolov-Miyaoka-Yau inequality) appears in
several places in the literature, notably in work of Barthel,
Hirzebruch and H\"ofer~\cite{bhh}. In fact, after proving
Theorems~\ref{thm:compact} and~\ref{thm:noncompact}, we realized that
the lattices $\Sc(p,\bar\sigma_4)$ with even values of $p$ ($p=4,6,8$)
appear explicitly on p.~215 of the book~\cite{bhh}. Moreover, that
table suggested to us to include the value $p=\infty$.

The Barthel-Hirzebruch-H\"ofer construction was later reinterpreted
and extended by Couwenberg, Heckman and Looijenga, and it turns out
that all lattices $\Sc(p,\bar\sigma_4)$, $p=3,4,5,6,8,12$ appear in
their paper (see Table~8.5 of~\cite{chl}). We refer to these groups as
CHL lattices of Shephard-Todd type $G_{24}$. The isomorphism of
Theorem~\ref{thm:iso} also implies the following.
\begin{thmintro}
  The CHL lattice of Shephard-Todd type $G_{24}$ for $p=4,5,6,8,12$
  are not arithmetic. They are pairwise incommensurable, and they are
  not commensurable to any Deligne-Mostow lattice.
\end{thmintro}
To the author's knowledge, the arithmetic structure of the CHL ball
quotients was never worked out in the literature, nor was their
commensurability relations with other known ball quotients. It seems to
have been known to Deligne and Mostow, at least conjecturally, see
page~181 in~\cite{delignemostowbook} (Deligne and Mostow also refer to
the work of Yoshida~\cite{yoshidaDE}).

Note that part of our proof (the existence of a complex hyperbolic
uniformization for the pairs in Table~\ref{tab:pairs}) could in fact
be avoided, if we were to simply quote the results of Couwenberg,
Heckman, Looijenga. However we believe our proof, which is fairly short
and elementary, gives an interesting different perspective on this
special case of the results in~\cite{chl}.

More of the lattices that appear in~\cite{thealgo} can be treated with
techniques similar to the ones in this paper, see for
instance~\cite{abelian}.\\

\textbf{Acknowledgements:} I would like to thank St\'ephane Druel for
his invaluable help with some of the technical details needed in the
proof, Philippe Eyssidieux for his interest in this work, and Michel
Brion for pointing to precious references about the invariants of the
automorphism group of the Klein quartic. I also thank the referees for
various suggestions that helped improve the manuscript.

\section{Description of the groups $\Sc(p,\bar\sigma_4)$} \label{sec:s4c}

In order to describe the group $\Sc(p,\bar\sigma_4)$, rather than
using the original description, we use the characterization of this
lattice given in Proposition~\ref{prop:char_s4c}. The
characterization follows immediately from the results in~\cite{dpp}
(see also~\cite{thealgo}), we briefly review the argument.

In the next statement, we take $R_1,J$ to be elements of
$SU(2,1)$. Suppose that $R_1$ is a complex reflection with eigenvalues
$u^2,\bar u,\bar u$, where $u=e^{2\pi i/3p}$ and $J$ is a regular elliptic
element of order 3. We write $R_2=JR_1J^{-1}$, $R_3=J^{-1}R_1J$,
and $G$ for the group generated by $R_1$ and $J$.
\begin{prop}\label{prop:char_s4c}
  Suppose $G$ is a lattice, $(R_jR_k)^2=(R_kR_j)^2$ for all $j\neq k$,
  and $R_1J$ has order 7. Then $G$ is conjugate to
  $\Gamma_p=\Sc(p,\bar\sigma_4)$ (and $p$ is equal to $3, 4, 5, 6, 8$, or
  $12$).
\end{prop}
\begin{pf}
  From the braid relation $(R_jR_k)^2=(R_kR_j)^2$, it follows that the
  eigenvalues of $R_jR_k$ are roots of unity, namely $\bar u^2, \pm i
  u$ (see for instance Proposition~2.3 of~\cite{thealgo}). Clearly the
  eigenvalues of $R_1J$ are roots of unity as well. 

  It follows that the pair of matrices $(R_1,J)$ can be simultaneously
  conjugated to standard generators of a complex hyperbolic sporadic
  group, i.e. up to complex conjugation and multiplication by a cube
  root of unity, $\tr(R_1J)$ can be assumed to be one of the 10 values
  $\sigma_j$, $j=1,\dots,10$ listed in~\cite{thealgo}.

  Inspection of the values of the order of $R_1J$ and the braid length
  $\textrm{br}(R_j,R_k)$ for these values of $\tr(R_1J)$ show that we may
  assume $\tr(R_1J)=\sigma_4$ or $\bar\sigma_4$. Now the group
  $\Sc(p,\sigma_4)$ is a subgroup of $PU(2,1)$ only for $p=4,5$ or
  $6$, and it is known to be non-discrete for each of these three
  values. The groups $\Sc(p,\bar\sigma_4)$ are subgroups of $PU(2,1)$
  for $p>2$, and if so they are not discrete for $p\neq
  3,4,5,6,8,12$. Both non-discreteness statements just mentioned
  follow from Theorem~9.1 of~\cite{dpp}.
\end{pf}

\section{Weighted projective plane} \label{sec:weighted}

The group $\Gamma_2=\Sc(2,\bar\sigma_4)$ was not studied
in~\cite{dpp}, because it does not act on the complex hyperbolic
plane. On the other hand, the matrices given there make sense for
$p=2$, and in that case the group generated by $R_1$ and $J$ preserves
a definite Hermitian form (which is unique up to scaling). In other
words, the group can be seen as a discrete (in fact finite) group of
isometries of the Fubini-Study metric on $\P^2$.

We omit the proof of the following Proposition, since it is not needed
anywhere in the paper (it is included only for motivation).
\begin{prop}\label{prop:iso168}
  The group $\Gamma_2$ is isomorphic to the projective automorphism
  group of the Klein quartic, which is also the 24-th Shephard-Todd
  group $G_{24}$.
\end{prop}
From this point on, we work only with the automorphism group of the
Klein quartic, and we write it as $G$ (and we write $\P G$
for its projectivization).

\begin{prop}\label{prop:g168}
  The group $\P G$ is the finite simple group of order 168. The
  quotient $\quotr{\P^2}{G}$ 
  is isomorphic as complex analytic
  orbifold to the pair $(X,\frac{1}{2}M)$ where $X=\P(2,3,7)$ and $M$
  is the image of the union of mirrors of reflections in
  $\P G$. Moreover,
  \begin{enumerate}[(i)]
    \item $M$ is an irreducible curve of degree 21 in
      $\P(2,3,7)$, with equation~\eqref{eq:discr} in a natural set of
      coordinates,
    \item $M$ contains precisely one singular point of $\P(2,3,7)$,
      which is its $A_1$ singularity,
    \item $M$ has two singular points $s_3$ and $s_4$ in the smooth
      part of $\P(2,3,7)$, where it has (analytic) local equations of
      the form $w_1^3=w_2^2$ and $w_1(w_1-w_2^2)=0$ respectively.
  \end{enumerate}
\end{prop}
We refer to the curve $M$ described in the proposition as the \emph{Klein
discriminant}.
Recall that $\P(2,3,7)$ is the quotient of $\C^3$ by the equivalence
relation $(z_1,z_2,z_3)\sim(\lambda^2 z_1,\lambda^3 z_2, \lambda^7
z_3)$ for all $\lambda\in\C^*$.  It has three singular points, given
in (weighted) homogeneous coordinates $z_1,z_2,z_3$ by the coordinate
axes, in other words $(z_1,z_2,z_3)=(1,0,0)$, $(0,1,0)$,
$(0,0,1)$. The corresponding singularities are of type $A_1$, $A_2$
and $\frac{1}{7}(2,3)$, respectively.

One way to describe the curve $M$ is to give an explicit equation,
which in a suitable set $(z_1,z_2,z_3)$ of coordinates can be taken to
be
\begin{align*}
z_3^3 + 27\cdot 64 z_2^7 
- 88z_1^2z_2z_3^2 + 16\cdot 63z_1z_2^4z_3 + 17\cdot 64z_1^4z_2^2z_3 \qquad \numberthis \label{eq:discr}\\
\qquad -256z_1^7z_3 - 128\cdot 469 z_1^3z_2^5
+ 43\cdot 512 z_1^6z_2^3 - 2048 z_1^9z_2=0.
\end{align*}
This curve may seem very mysterious, and its structure 
(irreducibility, nature of singular points) is not completely
obvious. To prepare for the proof of Proposition~\ref{prop:g168}, we
first briefly review some classical facts about the automorphism group
of the Klein quartic.

The Klein quartic is the complex curve given in homogeneous
coordinates $(x_1,x_2,x_3)$ for $\CP^2$ by $f(x_1,x_2,x_3)=0$ where
\begin{equation}\label{eq:kleinquartic}
f(x_1,x_2,x_3)=x_1^3x_2 + x_2^3x_3 + x_3^3x_1.
\end{equation}
It is well known that its automorphism group is the only simple
group of order 168, and that it can be written out explicitly as the
(projective) unitary group generated by the following three matrices:
$$
T=\left(\begin{matrix}
  \zeta&0&0\\
  0&\zeta^2&0\\
  0&0&\zeta^4
\end{matrix}\right)\quad 
J=\left(\begin{matrix}
  0&0&1\\
  1&0&0\\
  0&1&0
\end{matrix}\right)\quad
R=\left(\begin{matrix}
  a&b&c\\
  b&c&a\\
  c&a&b
\end{matrix}\right),
$$ 
where $\zeta=e^{2\pi i/7}$, $a=h(\zeta^4-\bar\zeta^4)$,
$b=h(\zeta^2-\bar\zeta^2)$, $c=h(\zeta-\bar\zeta)$, and
$h=i/\sqrt{7}=1/(\zeta^3+\zeta^5+\zeta^6-\bar\zeta^3-\bar\zeta^5-\bar\zeta^6)$.

\begin{rk}
  The group $\langle R,J,T\rangle$ is not quite generated by complex
  reflections, but it is a subgroup of index two in a reflection
  group, obtained from it by adjoining $-I$. The corresponding group
  of order 336 is the 24-th group in the Shephard-Todd
  list~\cite{shephardtodd}.
\end{rk}

Note that $R$ has order 2, and it is (the opposite of) a complex
reflection (it has eigenvalues $1$, $-1$, $-1$). We write $v_1$ for a
unit eigenvector with eigenvalue 1, which is orthogonal to the mirror
of $R$.
\begin{prop}\label{prop:finite}
  The matrices $T$ and $J$ generate a group $H$ of order 21. The orbit
  of the mirror of $R$ under this group has 21 elements, and the
  stabilizer of $\mathbb{C}v_1$ in $G$ has order 8. The group $G$ has
  order 168, and the conjugates of $R$ by elements of $H$ are the only
  reflections in $G$.
\end{prop}
\begin{pf}
  The first statement is obvious, since $J$ permutes the eigenspaces
  of $T$ cyclically. Let us denote by $V$ the orbit of $v_1$ under
  $H$. One easily checks that no element of $H$ fixes $v_1$, so $V$
  has 21 elements.

  Among these 21 vectors, one checks that 4 are orthogonal to $v_1$,
  and these come in two pairs $\{a_1,a_2\}$, $\{a_3,a_4\}$ of
  orthogonal vectors ($a_1$ and $a_2$ are orthogonal to each other,
  but they are not orthogonal to $a_3$ nor to $a_4$). Explicitly,
  $a_1=JTv_1$, $a_2=J^{-1}T^3v_1$, $a_3=J^{-1}T^4v_1$,
  $a_4=JT^{-1}v_1$.

  This implies that the stabilizer of $v_1$ in $G$ has order dividing
  8. One checks that $-TRTJRJ$ fixes $v_1$ and acts by switching $a_1$
  and $a_2$, and $-T^2JTRT$ fixes $v_1, a_1$ and changes $a_2$ to
  $-a_2$. Since these two transformations generate a group of order 8,
  the stabilizer of $v_1$ has order 8.

  This implies that $G$ has order $21\cdot 8=168$.
\end{pf}
\begin{pf} (of Proposition~\ref{prop:g168})
  Part~(i) follows from the explicit determination of generators for
  the ring of invariants for $G$. This was worked out in the
  nineteenth century by Klein~\cite{klein}, see
  also~\cite{weber},~\cite{springer}.

  It is not completely obvious that the
  equation~\eqref{eq:kleinquartic} of the Klein quartic is invariant
  under $R$, but it can be checked directly. This gives an invariant
  of degree 4. One then constructs another invariant from the Hessian
  $H(f)=\left( \frac{\partial^2 f}{\partial x_j\partial
    x_k}\right)_{j,k}$, namely $\Delta=\frac{1}{9}\det(H(f))$, which
  is homogeneous of degree 6, then
  $$
  C = \frac{1}{9}\det\left(
  \begin{matrix}
    H(f) & \gr(\Delta)\\
    \gr(\Delta)^t& 0
  \end{matrix}
  \right),
  $$
  which is homogeneous of degree 14, and finally
  $$
  K=\frac{1}{14}\det\left(\gr f,\gr\Delta,\gr C\right),
  $$
  which has degree 21.

  Note that the last polynomial is not invariant under $-I$ (but its
  square is). The ring of invariants for the group $\langle
  G,-I\rangle$, which is a reflection group, is generated by
  $f,\Delta,C$, so the quotient of $\CP^2$ is $\P(4,6,14)$ which is
  isomorphic to $\P(2,3,7)$.

  In particular, $K^2$ can be expressed as a polynomial in
  $f,\Delta,C$, which was also computed in the nineteenth century by
  Gordan~\cite{gordan}. This gives an equation for $D$, which is
  equation~\eqref{eq:discr}, in the coordinates
  $(z_1,z_2,z_3)=(f,\Delta,C)$ (see~\cite{springer} or p.529
  of~\cite{weber}).

  The irreducibility of $M$ follows from the fact that $G$ acts
  transitively on the set of mirrors of reflections in $G$ (see
  Proposition~\ref{prop:finite}). Part~(ii) is obvious from
  equation~\eqref{eq:discr}.

  The proof of part~(iii) relies on the detailed study of fixed points
  in $\CP^2$ of elements in $\P G$. One checks that there are 171
  fixed points of regular elements in the group, and these fixed
  points come in 5 orbits with various stabilizers, as listed in
  Table~\ref{tab:orbits} (see \S136 of~\cite{weber}). When the
  stabilizer is generated by complex reflections, we list its
  Shephard-Todd notation in the last row (in that case the
  corresponding orbit gives a smooth point of $\P(2,3,7)$). In
  particular, the group generated by $R$ and $J$ is generated by $R$
  and $JR$, and these two are indeed complex reflections. We denote by
  $M_1$ and $M_2$ the matrices that already appeared (up to their
  sign) in the proof of Proposition~\ref{prop:finite}, namely
  $M_1=TRTJRJ$ and $M_2=T^2JTRT$.
  \begin{table}[htbp]
    \begin{tabular}{c|ccccc}
      Notation               & $s_3$      & $s_4$       & $t_2$    & $t_3$  & $t_7$\\\hline
      \# mirrors             &  3         &  4          & 1        & 0      & 0\\
      \# stab                &  6         &  8          & 4        & 3      & 7\\
      \# orbit               &  28        &  21         & 42       & 56     & 24\\ 
      Gens                   &  $J,R$     & $M_1,M_2$   & $TR$     & $J$    & $T$\\
      ST                     & $G(3,3,2)$ & $G(2,1,2)$  & $\times$ & $\times$ & $\times$
    \end{tabular}
    \caption{The list of orbits of fixed points of regular elliptic
      elements in the group. For each point, we give the number of mirrors
      through that point, the order of and generators for each stabilizer
      in $G$, and if applicable the Shephard-Todd notation, if the
      stabilizer is generated by complex reflections.}\label{tab:orbits}
  \end{table}
 Note that the configuration of mirrors of reflections in $G$ has 28
 triple points, and 21 quadruple points, see page~96 of~\cite{bhh}.

 The local analytic structure of the singularities of $M$ can be
 understood by studying invariants for the $G(m,p,n)$ groups that
 appear here in Table~\ref{tab:orbits}. In suitable coordinates for
 $\C^2$, the group $G(3,3,2)$ is generated by
 $(z_1,z_2)\mapsto(z_2,z_1)$ and $(z_1,z_2)\mapsto(\omega z_1,
 \overline{\omega} z_2)$. The ring of invariants is the polynomial
 ring in $u_1=z_1z_2$, $u_2=z_1^3+z_2^3$, in other words the map
 $(z_1,z_2)\mapsto(u_1,u_2)=(z_1z_2,z_1^3+z_2^3)$ gives the quotient
 map $\C^2\rightarrow \C^2$ by the group action of $G(3,3,2)$.

 The mirrors of reflections in the group are $z_1=z_2$, $z_1=\omega
 z_2$, $z_1=\overline{\omega} z_2$, so the union of the mirrors has
 equation $z_1^3-z_2^3=0$, which has the same zero set as the
 invariant polynomial $(z_1^3-z_2^3)^2$. In terms of our invariant
 generators, the union of the mirrors is given by $u_2^2-4u_1^3=0$.
 This also gives the branch locus of the branched covering of the
 quotient map.

  Similarly, the group $G(2,1,2)$ is generated by
  $(z_1,z_2)\mapsto(z_2,z_1)$ and $(z_1,z_2)\mapsto(- z_1,z_2)$. The
  corresponding ring of invariants is the polynomial ring in
  $u_1=z_1^2z_2^2$ $u_2=z_1^2+z_2^2$. The mirrors of reflections in
  the group are $z_1=\pm z_2$, $z_1=0$, $z_2=0$, so the union of the
  mirrors has equation $z_1^2z_2^2(z_1-z_2)^2=0$, or equivalently
  $u_1(u_2^2-4u_1)=0$.
\end{pf}

One easily verifies that a given mirror of any reflection in the group
contains precisely 10 points that are fixed by elements in the group
other than the reflection itself and the central element. We deduce
the following result, which will be used later in the paper.
\begin{prop}\label{prop:chizero}
  Let $X=\P(2,3,7)$, and let
  $X_0=X\setminus(M\cup\{t_2,t_3,t_7\})$. Then $\chi(X_0)=0$.
\end{prop}
\begin{pf}
  This follows by computing the Euler characteristic of the complement
  $\tilde{X}_0$ in $\CP^2$ of the union of all fixed point sets of
  non-trivial elements in the group $\P G$, where the action of $\P G$
  is free. $X_0$ is obtained from $\CP^2$ by removing 21 copies of
  $\P^1\setminus\{10 pts\}$, as well as the 171 fixed points of
  regular elements (these include intersections of mirrors). This
  gives
  $$
  \chi(\CP^2) = 3 = \chi(\tilde{X}_0) + 21\cdot(2-10) + 171, 
  $$
  so $\chi(\tilde{X}_0)=0$, and $\chi(X_0)=\chi(\tilde{X}_0)/168=0$.
\end{pf}

In order to compute orbifold Chern numbers, we will need some basic
algebro-geometric properties of weighted projective space $\P(2,3,7)$,
all of which can be derived from the general theory of toric varieties
(see~\S3.4 and~\S4.3 of~\cite{fulton}, or~\cite{cls}).
\begin{prop}\label{prop:pabc}
  Let $X=\P(a_1,a_2,a_3)$ be a well-formed weighted projective space,
  i.e. such that every pair of integers in $\{a_1,a_2,a_3\}$ is
  relatively prime.
\begin{enumerate}
\item The Picard group ${\rm Pic}(X)$ is isomorphic to $\Z$, we denote
  by $H$ its positive generator.
\item The canonical divisor $K_X$ is equivalent to $-(a_1+a_2+a_3)H/a_1a_2a_3$.
\item $H^2=a_1a_2a_3$.
\end{enumerate}
\end{prop}

\section{Ball quotients via uniformization}

In order to construct orbifolds uniformized by the ball, we would like
to consider pairs $(X,(1-\frac{1}{p})M)$ for integers $p\geq 2$, where
$X=\P(2,3,7)$ and $M$ is the Klein discriminant (see
Proposition~\ref{prop:g168}). This idea is somewhat natural given the
families of lattices studied in~\cite{thealgo}, and also in view of
the work by Couwenberg, Heckman and Looijenga. Indeed, they
constructed 1-parameter families of complex hyperbolic structures on
$\P^2\setminus \mathcal{H}$, where the parameter corresponds to the
common rotation angle of the holonomy around all lines in
$\mathcal{H}$. Moreover, their structures actually descend to the
quotient by the action of the finite group $G_{24}$ (see section~6
of~\cite{chl}).

We have already seen that $(X,(1-\frac{1}{2})M)$ is an orbifold
uniformized by $\P^2$ (see Proposition~\ref{prop:g168}). Perhaps
surprisingly, we will see that the pair $(X,(1-\frac{1}{3})M)$ is an
orbifold as well, now uniformized by $\H^2$.  For higher values of
$p$, more drastic modifications of $X$ need to take place in order to
make the orbifold uniformizable. In order to explain these
modifications, we start by reviewing the logarithmic
Bogomolov-Miyaoka-Yau inequality.

\subsection{The logarithmic Bogomolov-Miyaoka-Yau inequality} \label{sec:logBMY}

We review some results about identifying ball quotients by computing
ratios of Chern classes. It is well known that a smooth compact
complex surface of general type $X$ is a ball quotient if and only if
$c_1^2(X)=3c_2(X)$. Note that the latter relation between Chern
numbers holds for $X=\CP^2$, so it also holds for compact ball
quotients by the Hirzebruch proportionality principle.

The fact that a surface of general type whose Chern numbers satisfy
this relation is indeed a ball quotient follows from the solution of
the Calabi conjecture by Aubin and Yau. Indeed, on a surface of
general type, there exists a K\"ahler-Einstein metric with negative
constant, and a classical computation shows that if
$c_1^2(X)=3c_2(X)$, then that metric actually has constant holomorphic
sectional curvature (and the constant can be taken to be equal to
$-1$). It then follows that $X$ can be written as $\quot{\H^2}{\Gamma}$
for some \emph{torsion-free} lattice $\Gamma$ in
${\rm Bihol}(\H^2)=\PU(2,1)$.

Ball quotients $\quot{\H^2}{\Gamma}$
where $\Gamma$ is a lattice
that is not torsion-free are in general not smooth surfaces, but they
are smooth as orbifolds. The corresponding map $\H^2\rightarrow
\quot{\H^2}{\Gamma}$ is then a branched cover, with branch locus
given by the set of fixed points of non-trivial elements in
$\Gamma$. The quotient is smooth near the orbit $\Gamma\cdot x$ if and
only if the stabilizer of $x$ in $\Gamma$ is generated by complex
reflections, by a theorem of Chevalley~\cite{chevalley}.

One naturally expects that the above characterization of ball
quotients should hold for orbifold logarithmic pairs $(X,D)$ (here $D$
plays is a $\Q$-divisor that corresponds the ramification divisor of
the quotient map from the ball), provided we replace the usual Chern
numbers by orbifold Chern numbers. A statement along those lines was
indeed proved by Kobayashi, Nakamura and Sakai~\cite{kns}, see
also~\cite{kobayashi},~\cite{bhh} and~\cite{tretkoff}.

The setting is as follows. Let $(X,D)$ be a pair with $X$ a compact
complex normal surface, and $D$ is a $\Q$-divisor of the form
$D=\sum_j(1-\frac{1}{b_j})D_j$, each $D_j$ being distinct irreducible
curves and $b_j=2,3,\dots$ or $\infty$ (in the latter case, we set
$1/\infty=0$).
We assume that the pair $(X,D)$ has at worst log-canonical
singularities, and that $K_X+D$ ample. Then the Kodaira dimension
$\kappa(X,D)$ is equal to two, the log-canonical ring
$R(X,D)=\oplus_{m\geq 0} H^0(X,m(K_X+D))$ is finitely generated, and
$X={\rm Proj}(R(X,D))$ is its own log-canonical model.

Let $X_0$ be obtained from $X$ by removing the components $D_j$ such
that $b_j=\infty$, as well as the singular points of $(X,D)$ that are
\emph{not} log-terminal, and let $D_0=D\cap X_0$. The the main result
prove by Kobayashi, Nakamura and Sakai in~\cite{kns} is the following.
\begin{thm}\label{thm:kns}
Under the above hypotheses, the pair $(X_0,D_0)$ is an orbifold,
$c_1(X_0,D_0)^2 \leq 3 c_2(X_0,D_0)$ and equality holds if and only if
$(X_0,D_0)$ is a quotient $\quot{\H^2}{\Gamma}$ for some lattice
$\Gamma\subset \PU(2,1)$.
\end{thm}
Note that $c_2(X_0,D_0)$ makes sense because the pair $(X_0,D_0)$ is
an orbifold, and it is often called the orbifold Euler
characteristic. The left hand side $c_1(X_0,D)$ can also be computed
as the self-intersection $(K_{X}+D)^2$ of the orbifold canonical
divisor (we will also refer to the latter as the log-canonical divisor
of the pair $(X,D)$).

We briefly recall the definition of log-canonical and log-terminal
singularities. Let $p$ be a singular point of $(X,D)$, and let
$\mu:\widetilde{X}\rightarrow X$ be a log-resolution of that
singularity. Here $\widetilde{D}$ denotes the strict transform of $D$,
and we write $\widetilde{E}$ for the exceptional set, and write
$E_\alpha$ for its irreducible components. We assume $K_X+D$ is
$\Q$-Cartier, and that the resolution is ``good'', i.e. that $G$ has
normal crossings, where the divisor $G$ obtained by writing
$$
K_{\widetilde{X}} + \widetilde{D} = \mu^*(K_X+D) + G 
$$ 
The divisor $G$ can be written as $\sum_\alpha g_\alpha E_\alpha$,
and then the singularity is called log-canonical (resp. log-terminal)
if $g_\alpha\geq -1$ for all $\alpha$ (resp. $g_\alpha>-1$ for all
$\alpha$, and $b_i<\infty$ for all $i$). This condition is known to
be independent of the good resolution $\mu$.

As mentioned above, in order to apply the Kobayashi-Nakamura-Sakai
result, there are three things to check, each handled in one of the
following sections. In section~\ref{sec:singularities}, we study the
singularities of the pairs. In section~\ref{sec:miyaoka} we verify
that equality holds in the orbifold Bogomolov-Miyaoka-Yau
inequality. In section~\ref{sec:ample} we check that the
log-canonical divisors of the relevant pairs are ample.

\subsection{Description of the logarithmic pairs}

The pairs $(X,(1-\frac{1}{p})M)$ do not usually satisfy the technical
hypotheses of the result Kobayashi, Nakamura and Sakai (see
section~\ref{sec:logBMY}). The pairs we will use to produce ball
quotients are obtained from $(X,(1-\frac{1}{p})M)$ by suitable
birational modifications.

Later, we will describe the relevant birational transformations
directly on the level of $\P(2,3,7)$, see
section~\ref{sec:singularities}, but for now we give an equivariant
description on the level of $\CP^2$, seen as a branched covering of
$\P(2,3,7)$ of degree 168. It is well known that the configuration of
mirrors of the group $G$ in $\CP^2$ has 49 singular points, that come
in 21 quadruple points (all in the same $G$-orbit), and 28 triple
points (also forming a $G$-orbit), see p.~211 in~\cite{bhh}, for
instance.
\begin{dfn}\label{dfn:xyz}
  Let $\widetilde{X}$ be $\CP^2$, and denote by $\widetilde{Y}$ the
  blow-up of $\widetilde{X}$ at every point of the $G$-orbit of
  quadruple points, and by $\widetilde{Z}$ the blow-up of
  $\widetilde{X}$ at all singular points of the union of mirrors in
  $G$ (all quadruple and triple points).  Let
  $X=\quotr{\widetilde{X}}{G}$, $Y=\quotr{\widetilde{Y}}{G}$ and
  $Z=\quotr{\widetilde{Z}}{G}$. We denote by $E$ the image in $Y$ (or
  $Z$) of any exceptional divisor in $\widetilde{Y}$ (or
  $\widetilde{Z}$) corresponding to quadruple mirror intersections.
  We denote by $F$ the image in $Z$ of any exceptional divisor in
  $\widetilde{Z}$ corresponding to triple mirror intersections. The
  divisor $E$ (resp. $F$) lies over a point of $\P(2,3,7)$ which we
  denote by $s_4$ (resp. $s_3$).  We denote by $M$ the image in $X$
  (or $Y$, or $Z$) of the union of mirrors in $\CP^2$.
\end{dfn}
Of course, $X$ is just $\P(2,3,7)$ (see
Proposition~\ref{prop:g168}). Note that $Y$ (resp. $Z$) is not the
usual blow-up of $X$ at $s_3$ (resp. at $s_3$ and $s_4$).  Indeed, the
image of the blown-up divisors are smooth points of $\P(2,3,7)$, but
$Y$ has a singular point in the image of the exceptional locus in
$\widetilde{Y}$.  Similarly, $Z$ has a singular point in the image of
each $G$-orbit of exceptional divisors in $\widetilde{Z}$.

\begin{rk}
  The singularities of $Y$ and $Z$ (see Definition~\ref{dfn:xyz}) can
  be studied by using explicit coordinates on the blow-up of $\C^2$ at
  a given point, and looking at the linearized action at every fixed
  point of the action of $G(2,1,2)$ (resp. of $G(3,3,2)$) on the
  blow-up. Since it is slightly cumbersome, we will bypass this
  verification and define explicit birational maps $Y\rightarrow X$
  and $Z\rightarrow X$, see section~\ref{sec:singularities}.
\end{rk}

In order to construct ball quotients, we consider the pairs given in
Table~\ref{tab:pairs}.
\begin{table}[htbp]
\begin{itemize}
  \item $X^{(3)}=X$, $D^{(3)}=(1-\frac{1}{3})M$.
  \item $X^{(4)}=X$, $D^{(4)}=(1-\frac{1}{4})M$.
  \item $X^{(5)}=Y$, $D^{(5)}=(1-\frac{1}{5})M+(1-\frac{1}{10})E$.
  \item $X^{(6)}=Y$, $D^{(6)}=(1-\frac{1}{6})M+(1-\frac{1}{6})E$.
  \item $X^{(8)}=Z$, $D^{(8)}=(1-\frac{1}{8})M+(1-\frac{1}{4})E+(1-\frac{1}{8})F$.
  \item $X^{(12)}=Z$, $D^{(12)}=(1-\frac{1}{12})M+(1-\frac{1}{3})E+(1-\frac{1}{4})F$.
  \item $X^{(\infty)}=Z$, $D^{(\infty)}=M+(1-\frac{1}{2})E+(1-\frac{1}{2})F$.
\end{itemize}
\caption{The above pairs will be shown to give orbifolds that are
  uniformized by the ball (possibly after removing log-canonical
  singularities that are not log-terminal).}\label{tab:pairs}
\end{table}
In section~\ref{sec:singularities}, we will see that these pairs have
at worst log-canonical singularities. As mentioned in
section~\ref{sec:logBMY}, it is also important to describe the locus
$X_0^{(p)}$ of points where each pair $(X^{(p)},D^{(p)})$ is actually
log-terminal, since it gives the open set uniformized by the ball in
Theorem~\ref{thm:kns}. We will prove the following.
\begin{prop}\label{prop:terminal_locus}
  The pairs $(X^{(p)},D^{(p)})$ are log-terminal for $p=3,5,8,12$. The
  log-terminal locus of $(X^{(4)},D^{(4)})$ is given by
  $X_0^{(4)}=X\setminus\{s_4\}$. The log-terminal locus of
    $(X^{(6)},D^{(6)})$ is given by $X_0^{(6)}=Y\setminus\{s_3\}$. The
      log-terminal locus of $(X^{(\infty)},D^{(\infty)})$ is
      $X_0^{(\infty)}=Z\setminus M$.
\end{prop}

\section{Ball quotients}

The main goal of this section is to prove Proposition~\ref{prop:verif}
which, by Theorem~\ref{thm:kns}, implies
Theorems~\ref{thm:compact} and~\ref{thm:noncompact}.
\begin{prop}\label{prop:verif}
Let $(V,D)=(X^{(p)},D^{(p)})$ be as in Table~\ref{tab:pairs} for
$p=3,4,5,6,8,12$ or $\infty$. Then
\begin{enumerate}
  \item The pair $(V,D)$ has at worst log-canonical singularities, and
    these are log-terminal singularities, except at $s_4$ when $p=4$,
    at $s_3$ when $p=6$, and along the divisor $M$ when $p=\infty$. 
  \item We have the equality $c_1(V,D)^2=3c_2(V_0,D_0)$, where
    $(V_0,D_0=D\cap V_0)$ is the complement in $(V,D)$ of the non
    log-terminal locus.
  \item The log-canonical divisor $K_V+D$ is ample, and the pair
    $(V,D)$ is its own log-canonical model.
\end{enumerate}
\end{prop} 
Among the items in Proposition~\ref{prop:verif}, the first item
guarantees that $(V_0,D_0)$ is an orbifold, thanks to the
classification of log-canonical singularities of pairs (we will also
describe the orbifold structure explicitly, since this is needed in
order to compute $c_2(V_0,D_0)$, see section~\ref{sec:c2}). The last
two items guarantee that the corresponding orbifold is uniformized by
the ball, by Theorem~\ref{thm:kns}.

Each part of Proposition~\ref{prop:verif} will be treated in a
separate section, see
sections~\ref{sec:singularities},~\ref{sec:miyaoka},~\ref{sec:ample}.

\subsection{Study of the singularities} \label{sec:singularities}

\subsubsection{Singularities of $X$} \label{sec:X}

Throughout this section, we write $\lambda=1-\frac{1}{p}$ (this is the
coefficient of $M$ in the orbifold canonical divisor). We first
consider the case $p=3$, and the singularities of $(X,D)$, where
$X=\P(2,3,7)$, and $D=(1-\frac{1}{3})M$. The singular points outside
$M$ are log-terminal, since they are quotient singularities (see
Proposition 4.18 in~\cite{kollarmori}, for instance).

There are three points to consider, namely $t_2$, which is the $A_1$
singularity of $\P(2,3,7)$, and $s_3$ and $s_4$, which are smooth
points of $X$. Near the $A_1$-singularity, there is nothing to show
since the pair $(X,\lambda M)$ is isomorphic to the quotient of $\C^2$
by a cyclic group generated by a diagonal map whose square is a
complex reflection, so this fits in the classification of
log-canonical singularities (see part~(1) of Theorem~5 in~\cite{kns}).

Near $s_3$ (resp. $s_4$), the curve $M$ has a local analytic equation
of the form $z_1^2=z_2^3$ (resp. $z_1(z_1-z_2^2)=0$). One way to check
the log-canonical character of these singularities is to identify
these curves as branch loci of suitable complex reflection groups
(see~\cite{shephardtodd} and~\cite{bannai}). We give a direct proof,
and along the way we derive formulas that will be used when computing
$c_1^2(X^{(p)},D^{(p)})$.

We consider a local resolution of the pair $(X,D)$ at $s_4$, which is
given by blowing up $s_4$, so that the strict transform of $M$ has two
transverse branches, and then blowing up that transverse intersection
once more. We denote by $\pi_4:X_4\rightarrow X$ the corresponding
composition of blow-ups, by $E_2$ the strict transform of the first
exceptional divisor (which is a $(-2)$-curve), and by $E_1$ the second
exceptional divisor (which is a $(-1)$-curve).

Since $K_{X_4}=\pi_4^*K_X+2E_1+E_2$ and
$\pi_4^*M=\widetilde{M}+4E_1+2E_2$ (where $\widetilde{M}$ denotes the
strict transform of $M$), we have
\begin{equation}\label{eq:kx4} 
K_{X_4}+\lambda \widetilde{M}=\pi_4^*(K_X+\lambda
M)+(2-4\lambda)E_1+(1-2\lambda)E_2,
\end{equation}
so the pair is log-canonical at $s_4$ if and only if
$\lambda\in[0,\frac{3}{4}]$, which is the case for $p\leq 4$. For
$p=4$, we get a log-canonical singularity which is not log-terminal.

We now consider a local resolution of the pair $(X,D)$ at $s_3$. In this
case, three successive blow-ups are needed; we blow-up the cusp of
$M$, the proper transform is then tangent to the exceptional
divisor. We then blow-up the point of tangency, which makes the
intersection transverse, and then blow-up that transverse
intersection. We denote by $\pi_3:X_3\rightarrow X$ the corresponding
sequence of blow-ups.

The exceptional set is a chain of three copies of $\P^1$, with
self-intersections $-2$, $-1$, $-3$, we denote the corresponding
curves by $F_2,F_1,F_3$ respectively. The corresponding formula is
\begin{equation}\label{eq:kx3}
K_{X_3}+\lambda \widetilde{M} =
  \pi_3^*(K_X+\lambda M)+(4-6\lambda)F_1+(2-3\lambda)F_2+(1-2\lambda)F_3,
\end{equation}
so we get a log-canonical singularity for $p\leq 6$ (which is not
log-terminal for $p=6$).

\subsubsection{Singularities of $Y$} \label{sec:Y}

The space $Y$ has a birational map $Y\rightarrow X$ which we now
describe. One simply does the same blow-up of $s_4$ as described in
section~\ref{sec:X}, and then contract the $(-2)$-curve. Note that
this contraction produces an $A_1$-singularity (this follows from the
uniqueness of the minimal resolution of surface singularities).

As before, we denote by $\pi:X_4\rightarrow X$ the composition of
three blow-ups described above, and we write $\gamma_4:X_4\rightarrow
Y$ for the contraction, $E=\gamma_4(E_1)$ and $\varphi_4:Y\rightarrow
X$ for the corresponding birational map, see the left diagram
in~\eqref{eq:diagram}.

\begin{equation}\label{eq:diagram}
\begin{gathered}
\xymatrix{
 **[l]\widetilde{M}\subset X_4  \ar[d]_{\gamma_4} \ar[dr]^{\pi_4} \\
 **[l] M'\subset Y \ar[r]_{\varphi_4} & **[r]X\supset M } 
\qquad 
\xymatrix{
X_3 \ar[d]_{\gamma_3} \ar[dr]^{\pi_3} \\
Z_3 \ar[r]_{\varphi_3} & X }
\end{gathered}
\end{equation}
By construction, $\gamma_4$ gives a good resolution of the pair
$(Y,D)$, where $D=\lambda M'+\mu E$. Here we denote by $\widetilde{M}$
the strict transform of $M$ in $X_4$, and by $M'$ the image of
$\widetilde{M}$ under $\gamma_4$. We also write $\lambda=1-\frac{1}{p}$ (resp.
$\mu=1-\frac{1}{m}$) for the coefficient of $\widetilde{M}$ (resp. $E$) in the
relevant log-canonical divisor, see Table~\ref{tab:pairs} for $p\geq
5$.

Note that $\gamma_4^*E=E_1+\frac{1}{2}E_2$. Indeed,
$\gamma_4^*E=E_1+aE_2$ for some $a\in \Q$, but $0=\gamma_4^*E\cdot
E_2=(E_1+aE_2)\cdot E_2=1-2a$, so $a=1/2$. In particular, we get
\begin{equation}\label{eq:e2}
E^2=E_1\cdot \gamma_4^*E=E_1\cdot(E_1+\frac{1}{2}E_2)=-1/2.
\end{equation}
We have $\gamma_4^*M'=\widetilde{M}+bE_2$ for some
$b\in\Q$. Intersecting both sides with $E_2$, we get
$\widetilde{M}\cdot E_2-2b=0$, but $\pi_4^*M=\widetilde{M}+4E_1+2E_2$,
so $\widetilde{M}\cdot E_2=0$, which gives $b=0$.  Finally, we note
that $K_{X_4}=\gamma_4^*K_Y$, since $A_1$ singularities are
crepant. For completeness, we mention that $K_{X_4}=\pi^*K_X+2E_1+E_2$,
and $K_Y=\varphi_4^*K_X+2E$.

Now we have
$$
K_{X_4}+\lambda \widetilde{M}+\mu E_1 = \gamma_4^*(K_Y+\lambda M'+\mu E) - \frac{\mu}{2} E_2,
$$ 
so the pair is log-terminal provided $\lambda,\mu\in [0,1[$,
    $-\frac{\mu}{2}>-1$ (only log-canonical if equality holds in some of
    these strict inequalities).

This shows that pairs the $(Y,D)$ corresponding to $p=5,6$ have
log-terminal singularities above $s_4$ (in fact this remains true for
the pairs $(Z,D)$ corresponding to $p=8,12$, that will be introduced
in the next section, since they have the same local structure near
$s_4$).

\subsubsection{Singularities of $Z$} \label{sec:Z}

In a similar way, the space $Z$ has a birational map
$\varphi_3:Z\rightarrow X$. Near $s_4$, we perform the same sequence
of blow-ups then contraction as in section~\ref{sec:Y}.

Near $s_3$, we perform three successive blow-ups as in
section~\ref{sec:X}, then contract the $(-2)$-curve (producing an
$A_1$ singularity) and $(-3)$-curve (producing a singularity of type
$\frac{1}{3}(1,1)$). Once again, the identification of the type of
singularities after contraction follows from the uniqueness of the
minimal resolution of surface singularities, and the knowledge of a
suitable model resolution. As a model, one can take the resolution of
the cone over the rational normal curve of degree $d=2$ or $3$, which
gives the Hirzebruch surface $\F_d$, the exceptional locus being a
$(-d)$-curve.

\begin{figure}[htbp]
  \epsfig{figure=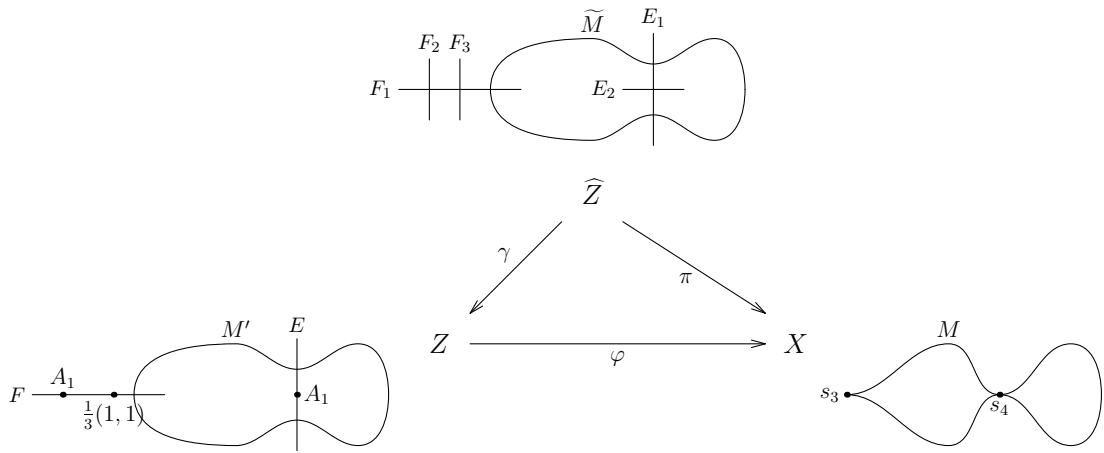,width=0.9\textwidth}
  \caption{Schematic picture of the space $Z$, which is obtained from
    $X$ by performing successive blow-ups at $s_3$ and $s_4$, then
    contracting the curves with self-intersection $(-2)$ or
    $(-3)$. This picture corresponds to the cases $p=8$ and $12$ (for
    $p=5$ or $6$, $s_3$ does not get blown-up).}\label{fig:klein}
\end{figure}

We denote by $\widehat{Z}$ the space $X$ blown-up twice at $s_4$ and
three times at $s_3$ (in the precise way that we just described), and
by $Z$ the space obtained from $\widehat{Z}$ by contracting the
exceptional curves with self-intersection $(-2)$ or $(-3)$. We denote
by $\pi:\widehat{Z}\rightarrow X$ the composition of blow-ups, by
$\gamma:\widehat{Z}\rightarrow Z$ the contraction, and by
$\varphi:Z\rightarrow X$ the corresponding morphism.  The
corresponding situation is illustrated in Figure~\ref{fig:klein}.
As mentioned in the end of the previous section, the pair $(Z,D)$ is
log-canonical above $s_4$, so we need only work locally (in the
analytic topology) around $s_3$. Hence we will work with $X_3$ rather
than $\widehat{Z}$.

We denote by $\pi_3:X_3\rightarrow X$ the sequence of blow-ups above
$s_3$, and we write $\gamma_3:X_3\rightarrow Z_3$ for the contraction,
$F=\gamma_3(F_1)$ and $\varphi_3:Z_3\rightarrow X$ for the
corresponding map (see the right part of
diagram~\eqref{eq:diagram}). As in section~\ref{sec:Y}, we denote by
$\widetilde{M}$ (resp. $M'$) the proper transform of $M$ in $X_3$
(resp. its push-forward in $Z_3$).  Again, $\gamma_3$ gives a good
resolution of the pair $(Z_3,D)$ near $s_3$, where $D=\lambda M'+\nu
F$. Here $\nu=1-\frac{1}{n}$ is the coefficient of $F$ in the relevant pair,
see Table~\ref{tab:pairs} for $p=8,12$.

We start by computing $F^2$, since it will be needed later when
computing $c_1^2$. We write $\gamma_3^*F=F_1+aF_2+bF_3$ for some $a,b\in
\Q$, and intersect both sides with $F_2$ or $F_3$, to get $1-2a$ and
$1-3b=0$. Then by the projection formula, 
\begin{equation}\label{eq:f2}
F^2=F_1\cdot\gamma_3^*F=-1+a+b=-\frac{1}{6}.
\end{equation}

The other relevant formulas are the following:
$$
\gamma_3^*F=F_1+\frac{1}{2}F_2+\frac{1}{3}F_3,
$$
$$
\pi_3^*M=\widetilde{M}+6F_1+3F_2+2F_3,
$$
and
$$
K_{X_3}=\pi_3^*K_X+4F_1+2F_2+F_3.
$$ 
Pushing the last formula forward by $\gamma_3$, we get
$K_{Z_3}=\varphi_3^*K_X+4F$, which gives
$$
K_{X_3}=\gamma_3^*K_{Z_3}-\frac{1}{3}F_3.
$$
One checks that $M'$ does not go through the singular points of $Z_3$ on $F$, so $\gamma_3^*M'=\widetilde{M}$ and
$$
K_{X_3}+\lambda \widetilde{M} + \nu F_1 = \gamma_3^*(K_{Z_3}+\lambda M'+\nu F) - \frac{\nu}{2} F_2 - \frac{1+\nu}{3} F_3.
$$ 
Hence the pair is log-terminal at points above $s_3$ if and only if
$\lambda,\nu\in[0,1[$, $\frac{\nu}{2}>-1$ and $-\frac{1+\nu}{3}>-1$,
    but only log-canonical if equality holds in some of these
    inequalities.

The result is that for $p>6$, all of the relevant pairs are
log-terminal above $s_3$.

\subsection{Equality holds in the Bogomolov-Miyaoka-Yau inequality} \label{sec:miyaoka}

\subsubsection{Computation of $c_1^2(X^{(p)},D^{(p)})$} \label{sec:c1}

For $p=3$ or $4$, by Proposition~\ref{prop:pabc}, we have
$K_X=-12(H/42)$, where $H$ is the positive generator of the Picard
group of $X$.
Since $M$ has weighted degree 21, $\lambda D=21\lambda (H/42)$, and
$$
(K_X+\lambda M)^2=\frac{1}{42}(-12+21\lambda)^2.
$$ 

For higher values of $p$, we push the formulas~\eqref{eq:kx4}
and~\eqref{eq:kx3} to $Y$ or $Z$, and we use the computation of $E^2$
and $F^2$ from section~\ref{sec:singularities} (see
equations~\eqref{eq:e2} and~\eqref{eq:f2}).

For $p=5$ or $6$, we use the map $\varphi_4:Y\rightarrow X$ from
section~\ref{sec:Y}. Equation~\eqref{eq:kx4} gives $ K_Y+\lambda
M'+\mu E = \varphi_4^*(K_X+\lambda M) + (2-4\lambda+\mu)E, $ hence
$$
(K_Y+\lambda M'+\mu E)^2 = \frac{1}{42}(-12+21\lambda)^2 + (-\frac{1}{2})(2-4\lambda+\mu)^2.
$$ 

For $p=8$, $12$ or $\infty$, we use the map $\varphi:Z\rightarrow
X$ from section~\ref{sec:Z}. Equations~\eqref{eq:kx4} and~\eqref{eq:kx3} give $ K_Z+\lambda M'+
\mu E + \nu F= \varphi^*(K_X+\lambda M) + (2-4\lambda+\mu)E +
(4-6\lambda+\nu)F, $ hence
$$
(K_Z+\lambda M'+\mu E+\nu F)^2 = \frac{1}{42}(-12+21\lambda)^2 + (-\frac{1}{2})(2-4\lambda+\mu)^2 + (-\frac{1}{6})(4-6\lambda+\nu)^2.
$$ 
We gather the corresponding numerical values, obtained from the
above formulas by taking specific values of $p,m,n$ corresponding to
Table~\ref{tab:pairs}. Recall that $\lambda=1-\frac{1}{p}$, $\mu=1-\frac{1}{m}$,
$\nu=1-\frac{1}{n}$.
\begin{table}[htbp]
  \begin{tabular}{| r | c c c c c c c |}
    \hline    
    $p$                      & 3 & 4 & 5  & 6 & 8 & 12 & $\infty$\\\hline
    $m$                      &   &   & 10 & 6 & 4 & 3  & 2\\\hline
    $n$                      &   &   &    &   & 8 & 4  & 2\\\hline
                             &   &   &    &   &   &    &   \\[-10pt]
    $c_1^2(X^{(p)},D^{(p)})$ & $\displaystyle\frac{2}{21}$ & $\displaystyle\frac{75}{224}$ & $\displaystyle\frac{141}{280}$ & $\displaystyle\frac{25}{42}$ & $\displaystyle\frac{297}{448}$ & $\displaystyle\frac{221}{336}$ & $\displaystyle\frac{3}{7}$\\[10pt]\hline
  \end{tabular}
\caption{Numerical values of $c_1^2(X^{(p)},D^{(p)})$, for various values of $p$.}\label{tab:c1}
\end{table}

\subsubsection{Computation of $c_2(X_0^{(p)},D_0^{(p)})$} \label{sec:c2}

Recall that we need to consider the pairs $(X^{(p)},D^{(p)})$ for
$p=3,5,8,12$, and $(X_0^{(p)},D_0^{(p)})$ for $p=4,6,\infty$. 

In order to compute $c_2(X_0^{(p)},D_0^{(p)})$, we compute the
orbifold Euler characteristic separately on each stratum with constant
isotropy group and sum the corresponding values. More precisely, we
write $X_0^{(p)}$ as a disjoint union $\sqcup_{s\in\mathfrak{S}}S$
where the $S$ have constant isotropy group, and compute
\begin{equation}\label{eq:eulerchar}
  \chi^{orb}(X_0^{(p)},D_0^{(p)})=\sum_{S\in\mathfrak{S}} \frac{\chi(S)}{|I_S|},
\end{equation}
where $\chi(S)$ is the usual topological Euler characteristic, and
$I_S$ is the isotropy group of an arbitrary point of $S$,
see~\cite{satake} (and
also~\cite{allendoerferweil},~\cite{mcmullengaussbonnet}).

Note that this formula is not exactly the same as the formula on the
right hand side of the inequality that appears in Theorem 12
of~\cite{kns}, which is closer to a Riemann-Hurwitz formula. One
reason why we use our formulation is that the Kobayashi-Nakamura-Sakai
formula is slightly ambiguous in our situation (the $d_i$ that is
stated in~\cite{kns} to be a number of singularities of a certain type
should be counted with multiplicity).

We break up the sum in equation~\eqref{eq:eulerchar} into summands
corresponding to each complex dimension $k=0,1,2$, and write
$$
\chi^{orb}=\chi^{orb}_0+\chi^{orb}_1+\chi^{orb}_2.
$$ 
By Proposition~\ref{prop:chizero}, the stratum with trivial
isotropy group has Euler characteristic 0, hence for every value of
$p$ we have
$$
\chi^{orb}_2(X_0^{(p)},D_0^{(p)})=0.
$$

For $\chi^{orb}_0$ and $\chi^{orb}_1$, we will get different formulas
depending on $p$. In fact the general form of our formulas depends on
whether $X^{(p)}$ is equal to $X$, $Y$ or $Z$ (see
section~\ref{sec:singularities}). We suggest the reader to keep
Figure~\ref{fig:klein} in mind, since it helps keep track of the
combinatorics and topology of the strata with constant isotropy
groups.

The most complicated istropy groups are the isotropy groups at the two
singular points $s_3$ and $s_4$ of the Klein discriminant. Since
$\P(2,3,7)$ is smooth at those points, the corresponding isotropy
groups must be 2-dimensional finite groups generated by complex
reflections, which were tabulated by Shephard-Todd,
see~\cite{shephardtodd}.

For each such singular point $s_k$ ($k=3$ or $4$), we consider a
``small'' contractible neighborhood $U_k$ of $s_k$ in $X$.
\begin{prop}\label{prop:braiding}
  For $k=3$ or $4$, and for $U_k$ small enough, we have
  $\pi_1(U_k\setminus M)=\langle
  \,a,b\,|\,(ab)^{k/2}=(ba)^{k/2}\,\rangle$.
\end{prop}
Recall that, for an odd integer $n$, $(xy)^{n/2}$ stands for an
alternating product $xyx\cdots yx$ with $n$ factors, and we call the
relation $(xy)^{n/2}=(yx)^{n/2}$ a braid relation of length $n$. For
short, when $x$ and $y$ braid with length $n$, we write
$\textrm{br}_n(x,y)$.\\
\begin{pf}
  This can be seen from standard arguments using projection onto one
  of the axes in $\C^2$ and studying the monodromy. Alternatively, one
  can use the description of the finite reflection stabilizers in the
  automorphism group of the Klein quartic, see the proof of
  Proposition~\ref{prop:g168} and Table~\ref{tab:orbits}. The
  corresponding local fundamental groups are tabulated
  in~\cite{bannai}.
\end{pf}
The local orbifold fundamental group at $s_3$, $s_4$ for the orbifold
$(X^{(p)},D^{(p)})$ are obtained from the groups in
proposition~\ref{prop:braiding} by adding the relations
$a^p=b^p=id$. Denote by $I_n(p)$ the group
$$
I_n(p)=\langle\, a,b\, |\, a^p, b^p, \textrm{br}_n(a,b)\, \rangle.
$$
\begin{prop}\label{prop:st}
  The groups $I_3(3)$, $I_3(4)$, $I_3(5)$, $I_4(3)$ are finite, of
  orders given by $|I_3(3)|=24$, $|I_3(4)|=96$, $|I_3(5)|=600$,
  $|I_4(3)|=72$. Each of these four $I_n(p)$ admits a faithful
  representation in $U(2)$ such that $a$ and $b$ are represented by
  complex reflections of order $p$. In particular, they are isomorphic
  to specific 2-dimensional primitive Shephard-Todd groups, namely
  $I_3(3)=G_4$, $I_3(4)=G_8$, $I_3(5)=G_{16}$ and $I_4(3)=G_5$.
\end{prop}
\begin{pf}
  This is explained in section~2.2 of~\cite{mostowpacific} for instance, see
  also the computation of the local fundamental group of the branch
  locus for 2-dimensional Shephard-Todd groups given in~\cite{bannai}.
\end{pf}

\noindent{\bf The cases $p=3,4$}

For $p=3$ or $4$, there is a unique stratum of dimension 1, namely
$M_0=M\setminus\{t_2,s_3,s_4\}$, where the isotropy group is cyclic of
order $p$. This is a $\P^1$ with four points removed (note that $s_4$
is a double point of $M$, see Figure~\ref{fig:klein}), hence
$\chi(M_0)=2-4=-2$ and
$$
\chi^{orb}_1=-\frac{2}{p}.
$$
In order to compute $\chi_0^{orb}$, we describe the local structure of
the corresponding orbifolds near points with isolated isotropy type
(i.e. points such that any neighboring point has a different isotropy
group).

\textsc{Local structure near $s_4$:} The local analytic structure of
the pair $(X,\frac{1}{3}M)$ near $s_4$ is given by the quotient of
$\C^2$ by the Shephard-Todd group $G_5$, see
Proposition~\ref{prop:st}, the corresponding isotropy group has
order 72.

For $p=4$, recall from section~\ref{sec:singularities} that the
singularity of the pair $(X^{(4)},D^{(4)})$ at $s_4$ is only
log-canonical, not log-terminal. So only
$X^{(4)}_0=X^{(4)}\setminus\{s_4\}$ carries an orbifold structure.

\textsc{Local structure near $s_3$:} The local analytic structure of
the pair $(X^{(p)},D^{(p)})$ near $s_3$ is given by the quotient of
$\C^2$ by the Shephard-Todd group $G_4$ for $p=3$, and by $G_8$ for
$p=4$ (see Poposition~\ref{prop:st}). These give isotropy groups of
order 24 and 96, respectively.

\textsc{Local structure near $t_2$:} Near the $A_1$-singularity, the
local model given by a regular elliptic element whose square is a
complex reflection of order $p$, hence we get isotropy group of order
$2p$.

\textsc{Local structure near $t_3,t_7$:} These two points are
not on the Klein discriminant curve, so we keep the orbifold structure
given by the quotient $\P(2,3,7)=\quotr{\P^2}{G_{24}}$.

The statements in the previous paragraphs are summarized in
Table~\ref{tab:chi0-34}.
\begin{table}
\begin{tabular}{|c|c|c|c|c|c|}
\hline
$p$ &   $t_2$   &  $t_3$  & $t_7$  &    $s_3$        & $s_4$          \\
\hline
2   & $\Z_{4}$  & $\Z_3$  & $\Z_7$ & $|G(3,3,2)|=6$  & $|G(2,1,2)|=8$ \\
3   & $\Z_{6}$  & $\Z_3$  & $\Z_7$ & $|G_4|=24$      & $|G_5|=72$     \\
4   & $\Z_{8}$  & $\Z_3$  & $\Z_7$ & $|G_8|=96$      & $\infty$       \\
\hline
\end{tabular}
\caption{Order of isotropy groups contributing to $\chi^{orb}_0$, for
  $p=2,3,4$.}\label{tab:chi0-34}
\end{table}
We get
$$
\chi^{orb}_0(X^{(3)},D^{(3)})=\frac{1}{2\cdot 3}+\frac{1}{24}+\frac{1}{72}+\frac{1}{3}+\frac{1}{7},
$$
$$
\chi^{orb}_1(X^{(3)},D^{(3)})=\frac{-2}{3},
$$
and finally
$$ 
c_2(X^{(3)},D^{(3)})=\chi^{orb}_0(X^{(3)},D^{(3)})+\chi_1^{orb}(X^{(3)},D^{(3)})=\frac{2}{63}.
$$ 
Similarly, 
$$ 
c_2(X^{(4)}_0,D^{(4)}_0)=(\frac{1}{2\cdot 4}+\frac{1}{96}+\frac{1}{3}+\frac{1}{7})+(\frac{-2}{4})=\frac{25}{224}.
$$ 
Comparing with the formulas in section~\ref{sec:c1}, we see that
$$
c_1^2(X_0^{(p)},D_0^{(p)})=3c_2(X_0^{(p)},D_0^{(p)})
$$
for both cases $p=3,4$.\\

\noindent{\bf The cases $p=5,6$}

Here and in the remainder of section~\ref{sec:c2}, with a slight abuse
of notation, we use the same notation for $M$ in $X$ and its proper
transform in $Y$ (or $Z$).

For $p=5,6$, there are two 1-dimensional strata, consisting of generic
points of $M$ and $E$, respectively. We denote by $M_0$ and $E_0$ the
corresponding non-compact curves. Just as in the cases $p=3,4$ we have
$\chi(M_0)=-2$.

Recall from section~\ref{sec:Y} that $Y$ is obtained from $X$ by
blowing up the point $s_4$ twice, then contracting the $(-2)$-curve in
the exceptional locus. This produces an $A_1$ singularity on $E$,
which we denote by $e_0$, and the divisor $E$ has two transverse
intersection points with $M$, which we denote by $e_1$ and $e_2$. Now
$E_0$ is a copy of $\P^1$ with three points removed (for a picture of
this, contract the exceptional divisor $F$ in Figure~\ref{fig:klein}).

In other words, we get
\begin{equation}\label{eq:chi1-56}
  \chi_1^{orb}(X_0^{(p)},D_0^{(p)})=\frac{-2}{p}+\frac{-1}{m},
\end{equation}
where $m=10$ (resp. $m=6$) if $p=5$ (resp. $p=6$).

In order to compute $\chi_0^{orb}$, we describe local models near each
point with special isotropy. 

\textsc{Local structure near $s_3$:} For $p=5$, the local model is
given the group generated by two complex reflections $a,b$ of order 5
that satisfy the braid relation $aba=bab$. This is the Shephard-Todd
group $G_{16}$, which has order 600.  For $p=6$, the corresponding
group would be infinite. In fact we have seen in
section~\ref{sec:singularities} that the pair $(X^{(6)},D^{(6)})$ is
not log-terminal at $s_3$, so the corresponding pair is not an
orbifold, and we need to remove that point in order to get a ball
quotient. In other words, $X^{(6)}_0=X^{(6)}\setminus\{s_3\}$,
$D^{(6)}_0=D^{(6)}\cap X^{(6)}_0$.

\textsc{Local structure near $e_1,e_2$:} These have abelian isotropy
groups generated by two complex reflections of order $p$ and $m$,
where $m=10$ (resp. $m=6$) if $p=5$ (resp. $p=6$), see
Table~\ref{tab:pairs}.

\textsc{Local structure near $e_0$:} Recall that this is an $A_1$
singularity on $E$. The isotropy group at this point is a regular
elliptic element whose square has order $m$, hence it has order $2m$.

\textsc{Local structure near $t_2$:} The isotropy group at this point
is the same as for $p=3,4$ (regular elliptic element whose square is a
complex reflection of order $p$), it has order $2p$.

\textsc{Local structure near $t_3,t_7$:} The order of these groups do
not change with $p$, they have order 3, 7
respectively.

We summarize the result of the previous paragraphs in
Table~\ref{tab:chi0-56}.
\begin{table}
\begin{tabular}{|c|c|c|c|c|c|c|}
\hline
$p$ &   $t_2$    &  $t_3$  & $t_7$  &    $s_3$        & $e_0$    & $e_1$, $e_2$  \\
\hline
5   & $\Z_{10}$  & $\Z_3$  & $\Z_7$ & $|G_{16}|=600$  & $2m=20$  &  $mp=50$      \\
6   & $\Z_{12}$  & $\Z_3$  & $\Z_7$ & $\infty$        & $2m=12$  &  $mp=36$      \\
\hline
\end{tabular}
\caption{Order of isotropy groups contributing to $\chi^{orb}_0$, for
  $p=5,6$.}\label{tab:chi0-56}
\end{table}
We then have
$$
\chi_0^{orb}(X^{(5)},D^{(5)})=\frac{1}{2\cdot 5}+\frac{1}{3}+\frac{1}{7}+
   \frac{1}{m}(\frac{2}{5}+\frac{1}{2})+\frac{1}{600},
$$
and
$$
\chi_0^{orb}(X^{(5)},D^{(5)})=\frac{1}{2\cdot 6}+\frac{1}{3}+\frac{1}{7}+
   \frac{1}{m}(\frac{2}{6}+\frac{1}{2}).
$$
Combining this with equation~\eqref{eq:chi1-56}, we get
$$ 
c_2(X^{(5)},D^{(5)})=\frac{47}{280},
$$
$$ 
c_2(X^{(6)}_0,D^{(6)}_0)=\frac{25}{126}.
$$ 
Comparing with the computations in section~\ref{sec:c1}, we see that
$c_1^2=3c_2$ for both $p=5,6$.\\

\noindent\textbf{The cases $p=8,12$}

Recall that to go from $X$ to $Z$, we do the same blow up and
contraction above $s_4$ as for $p<8$, and also a sequence of three
blow-ups above $s_3$, then contract a $(-2)$ and a $(-3)$-curves.
This produces two singular points, an $A_1$-singularity and a
singularity of type $\frac{1}{3}(1,1)$. We denote by $E$ the $\P^1$
that goes through the singular point corresponding to $s_4$, and $F$
the $\P^1$ that goes through the two singular points corresponding to
$s_3$ (see Figure~\ref{fig:klein}).

We use the same notation $e_0,e_1,e_2$ for special points on $E$ as in
the cases $p=5,6$, and we denote by $f_0$ the $A_1$ singularity on
$f$, by $f_1$ the $\frac{1}{3}(1,1)$ singularity, and by $f_2$ the
intersection point $F\cap M$.

Since the set of $E_0$ (resp. $F_0$) of generic points of $E$
(resp. $F$) is a $\P^1$ with three points removed, we have
$\chi(E_0)=\chi(F_0)=-1$. The generic isotropy groups on $M$ are
cyclic of $p$, cyclic of order $m$ on $E$ and cyclic of order $n$ on
$F$, where $(m,n)=(4,8)$ for $p=8$, $(m,n)=(3,4)$ for $p=12$ (see
Table~\ref{tab:pairs}). Hence we have
$$
\chi_1^{orb}=\frac{\chi(E_0)}{m}+\frac{\chi(F_0)}{n}+\frac{\chi(M_0)}{p}=
\frac{-1}{m}+\frac{-1}{n}+\frac{-2}{p}.
$$ 
The isotropy groups with special isotropy, corresponding to the
0-dimensional strata are listed in Table~\ref{tab:chi0-812}.
\begin{table}
\begin{tabular}{|c|c|c|c|c|c|c|c|c|}
\hline
$p$ &   $t_2$    &  $t_3$  & $t_7$  & $f_0$  & $f_1$   & $f_2$   & $e_0$   & $e_1$, $e_2$  \\
\hline
8   & $\Z_{16}$  & $\Z_3$  & $\Z_7$ & $2n=16$ & $3n=24$ & $np=64$ & $2m=8$ & $mp=24$   \\
12  & $\Z_{24}$  & $\Z_3$  & $\Z_7$ & $2n=8$ & $3n=12$ & $np=48$ & $2m=6$ & $mp=36$   \\
\hline
\end{tabular}
\caption{Order of isotropy groups contributing to $\chi^{orb}_0$, for
  $p=8,12$.}\label{tab:chi0-812}
\end{table}
We then get
$$
\chi_0^{orb}=\frac{1}{2\cdot p}+\frac{1}{3}+\frac{1}{7}
+\frac{1}{m}(\frac{2}{p}+\frac{1}{2})
+\frac{1}{n}(\frac{1}{p}+\frac{1}{2}+\frac{1}{3}),
$$
which gives
$$ 
c_2(X^{(8)},D^{(8)})=\frac{99}{448},\quad
c_2(X^{(12)},D^{(12)})=\frac{221}{1008}.
$$

\noindent\textbf{The case $p=\infty$}

In that case, we remove from $Z=X^{(\infty)}$ the curve $M$ above
$M$, and write $X^{(\infty)}_0=X^{(\infty)}\setminus M$,
$D^{(\infty)}_0=D^{(\infty)}\cap X^{(\infty)}_0$. In this case the
formulas are the same as for the cases $p=8,12$, but we set $p=\infty$
and $m=n=2$ (see Table~\ref{tab:pairs}). This gives
$$
\chi_1^{orb}(X^{(\infty)}_0,D^{(\infty)}_0)=\frac{-1}{m}+\frac{-1}{n},
$$
and
$$
\chi_0^{orb}(X^{(\infty)}_0,D^{(\infty)}_0)=\frac{1}{3}+\frac{1}{7}
+\frac{1}{m}(\frac{1}{2})
+\frac{1}{n}(\frac{1}{2}+\frac{1}{3}),
$$
and we get
$$
c_2(X^{(\infty)}_0,D^{(\infty)}_0)=\frac{1}{7}.
$$ 
Once again, this shows $c_1^2=3c_2$.

\subsection{Ampleness} \label{sec:ample}

The goal of this section is to prove that the pairs
$(X^{(p)},D^{(p)})$ as above have log-general type, i.e. that
$K_{X^{(p)}}+D^{(p)}$ is ample.

Note that this (would be false for $p=2$ and it) is easy for $p=3$
or $4$. Indeed, in that case, since $M$ has degree $21$, 
$$
K_X+\lambda M=(-12+21\lambda)H/42,
$$ 
where $H$ is the positive generator of ${\rm Pic}(X)$. The
coefficient $-12+21\cdot (1-\frac{1}{p})$ is $>0$ for $p\geq 3$
(whereas $-12+21\cdot (1-\frac{1}{2})=-\frac{3}{2}<0$). For $p\geq 5$
the question is slightly more subtle, since $X^{(p)}$ is not simply
given by $X$.

We first assume $p\geq 8$. Consider the map $\varphi:Z\rightarrow X$
constructed in section~\ref{sec:Z}.  Recall that
\begin{equation}\label{eq:klog}
K_Z + \lambda M' + \mu E + \nu F = \varphi^*(K_X+\lambda M) + (4-6\lambda+\mu) E + (2-4\lambda+\nu) F.
\end{equation} 
Note that the coefficients of $E$ and $F$ in the right hand side are
negative. 

We have $\varphi^*M=M'+6E+4F$, so $K_Z + \lambda
M' + \mu E + \nu F$ is ample if and only if
$\frac{21\lambda-12}{21}M'+A E+ BF$ is ample, where
$$
A=6\frac{21\lambda-12}{21}+(4-6\lambda+\mu)
$$ 
and 
$$
B=4\frac{21\lambda-12}{21}+(2-4\lambda+\nu).
$$ 

For $p=5$ or $6$, we do the same with $\varphi_4:Y\rightarrow
X$, but now
\begin{equation}\label{eq:klog2}
K_Y + \lambda M' + \mu E  = \varphi^*(K_X+\lambda M) + (4-6\lambda+\mu),
\end{equation} 
and we get a similar expression without $F$ (or in other words $B=0$).

The values of $A,B$ are listed in Table~\ref{tab:abvalues} (last three
columns), note in particular that $A,B>0$.
\begin{table}[htbp]
  \begin{tabular}{|r | ccccc |}\hline
    $p$ & 5 & 6 & 8 & 12 & $\infty$\\\hline
        &                     &                    &                    &                    & \\[-10pt]
    $A$ & $\di\frac{103}{70}$ & $\di\frac{59}{42}$ & $\di\frac{37}{28}$ & $\di\frac{26}{21}$ & $\di\frac{15}{4}$\\[2pt]
        &                     &                    &                    &                    &                  \\[-10pt]
    $B$ &                     &                    & $\di\frac{33}{56}$ & $\di\frac{13}{28}$ & $\di\frac{3}{14}$\\[8pt]
\hline
  \end{tabular}
\caption{Values of $A$ and $B$ for relevant values of $p$.}\label{tab:abvalues}
\end{table}

We now prove
\begin{prop}\label{prop:klognef}
  The divisor $K_{X^{(p)}}+D^{(p)}$ is ample for $p=3,4,5,6,8,12$ and
  $\infty$.
\end{prop}

\begin{pf}
  By the Nakai-Moishezon criterion, it is enough to check that the
  intersection $(K^{(p)}+D^{(p)})\cdot W$ is $>0$ for every irreducible curve $W$.

  Since $(K^{(p)}+D^{(p)})$ is numerically equivalent to
  $\frac{21\lambda-12}{21}M'+A E+ BF$ (or
  $\frac{21\lambda-12}{21}M'+A E$ if $p=5,6$), and the
  coefficients in this expression are all positive, it is enough to
  check that the intersection with $M'$, $E$ and $F$ are all
  $>0$.

  Indeed, if $W$ is distinct from $M'$, $E$ and $F$, its intersection
  with $(K^{(p)}+D^{(p)})$ is $\geq 0$. If it is 0, then $W$ projects
  to a curve in $X$ disjoint from $M$, but this is impossible since
  $X=\P(2,3,7)$ has Picard number one (see
  Proposition~\ref{prop:pabc}).

  For the intersection with $E$ or $F$, we use
  equation~\eqref{eq:klog} to get $(K_S + \lambda M' + \mu
  E)\cdot E = (4-6\lambda+\mu) E^2$ and $(K_S + \lambda M'
  + \mu E)\cdot F = (2-4\lambda+\nu) F^2$ (where $S$ is either $Y$ or
  $Z$, depending on the value of $p$). These numbers are both
  positive, since $E^2=-\frac{1}{2}<0$, $F^2=-\frac{1}{6}<0$ (see
  equations~\eqref{eq:e2} and~\eqref{eq:f2}) and for relevant values
  of $p$, we get $(2-4\lambda+\nu)<0$ and $(4-6\lambda+\mu)<0$, as
  shown in Table~\ref{tab:nef}.
  \begin{table}
    \begin{tabular}{|r|ccccc|}
      \hline
      $p$              & 5 & 6 & 8 & 12 & $\infty$ \\\hline
                       &   &   &   &    &  \\[-10pt]
      $2-4\lambda+\mu$ & $-\dfrac{3}{10}$ & $-\dfrac{1}{2}$ & $-\dfrac{3}{4}$ & $-1$ & $-\dfrac{3}{2}$\\[10pt]
      $4-6\lambda+\nu$ &           &           & $-\dfrac{3}{8}$ & $-\dfrac{3}{4}$ & $-2$\\[10pt]\hline 
    \end{tabular}
    \caption{Values of the coefficients of $E$ and $F$, which show that
      the intersection of the log-canonical divisor with $E$ and $F$ are
      positive.} \label{tab:nef}
  \end{table}

  Using $E^2=-\frac{1}{2}$ and $F^2=-\frac{1}{6}$ once again, we have
  that the intersection with $M'$, which is the same as the
  intersection with $\varphi^*M-6E-4F$ gives
  $$
  (21\lambda-12)\cdot 21\cdot\frac{1}{42} + (4-6\lambda+\mu) + 2(2-4\lambda+\nu),
  $$ 
  which is positive for $p=8,12$ and $\infty$ (the respective values
  are $23/16$, $23/24$ and $1/3$).
  
  For $p=5,6$, there is only one exceptional divisor $E$, so one simply
  needs to check $(2-4\lambda+\nu)<0$ (see the first two columns of
  Table~\ref{tab:nef}) and 
  $$
  (21\lambda-12)\cdot 21\cdot\frac{1}{42} + (4-6\lambda+\mu)>0,
  $$ 
  and this does indeed hold for $p=5,6$ (one gets $37/20$ and
  $17/12$, respectively).
\end{pf}

\section{Identification of the lattices} \label{sec:isomorphism}

In this section we prove Theorem~\ref{thm:iso}. The main ingredient is
a computation due to Naruki~\cite{naruki}, which describes the
fundamental group of the complement in $\P^2$ of the union of mirrors
of the reflections in the automorphism group of the Klein
quartic. Naruki actually computes the fundamental group of the
complement of the Klein discriminant in $\P(2,3,7)$, which is given as
follows (see Proposition~3.6 of~\cite{naruki}).
\begin{prop}
  The fundamental group of $X_0=X\setminus(M\cup\{t_2,t_3,t_7\})$ has
  a presentation of the form
  \begin{equation}\label{eq:pres_naruki}
  \langle\, \alpha,\delta\, |\, (\alpha\delta)^7,
  \textrm{br}_3(\alpha,\delta^2),
  \textrm{br}_4(\alpha,\delta^{-1}\alpha\delta)\,\rangle.
  \end{equation}
\end{prop}
Here we use the notation from~\cite{thealgo} and write
$\textrm{br}_n(a,b)$ for the relation $(ab)^{n/2}=(ba)^{n/2}$. For odd
$n$, $(ab)^{n/2}$ stands for a product of the form $aba\dots ba$ with
$n$ factors. In particular, $\textrm{br}_3(a,b)$ stands for $aba=bab$,
and this implies that $a$ and $b$ are conjugate, since
$a=(ba)b(ba)^{-1}$.

Following the proof given by Naruki, one easily sees that $\alpha$
corresponds to a loops that winds once around the Klein discriminant.

The inclusion $\iota:X_0\rightarrow X$ induces a surjective homomorphism
$\iota_*$ on the level of orbifold fundamental groups, such that
$\iota_*(\alpha)$ is a complex reflection of order $p$. It follows
that $\iota_*(\delta^2)$ is also a reflection of order $p$ (note that
$\alpha$ and $\delta^2$ are conjugate, since
$\textrm{br}_3(a,\delta^2)$).

In the remainder of this section, for the sake of readability, we
abuse notation and simply write $\gamma$ for $\iota_*(\gamma)$. We
also write $\bar\alpha$, $\bar\delta$ for $\alpha^{-1}$,
$\delta^{-1}$, respectively.

We will chose an isomorphism such that $\delta$ maps to the element
$S_1$ mentioned in equation~(10) of~\cite{dpp2} (which is a squareroot
of $R_1$), and $\alpha$ to $R_3^{-1}R_2R_3$.
\begin{prop}\label{prop:key}
  The elements $J=(\alpha\delta)^2\alpha\bar\delta$ and
  $R_1=\delta^2$ generate a group isomorphic to
  $\Sc(p,\overline{\sigma}_4)$.
\end{prop}
\begin{pf}
  One checks using the presentation~\eqref{eq:pres_naruki} that $J^3$
  has order 3. Indeed, 
  \begin{eqnarray*}
  &((\alpha\delta)^2\alpha\bar\delta)^2
    =\alpha\delta^2(\bar\delta\alpha\delta\alpha)^2\delta\alpha\bar\delta
    =\alpha\delta^2(\alpha\bar\delta\alpha\delta)^2\delta\alpha\bar\delta
    =(\alpha\delta^2\alpha)\bar\delta\alpha\delta\alpha\bar\delta(\alpha\delta^2\alpha)\bar\delta\\
  &=(\delta^2\alpha\delta^2)\bar\delta\alpha\delta\alpha\bar\delta(\delta^2\alpha\delta^2)\bar\delta
    =\delta^2(\alpha\delta)^4
    =\delta^2(\bar\delta\bar\alpha)^3
    =((\alpha\delta)^2\alpha\bar\delta)^{-1}.
  \end{eqnarray*}

  In the orbifold fundamental group seen as a subgroup of $PU(2,1)$,
  this element must be a regular elliptic element (this follows from
  the action of the group $G_{24}$ on $\P^2$, recall that every
  complex reflection in that group has order 2).  Since $\alpha\delta$
  has order 7, so does $P=R_1J=\delta(\delta\alpha)^3\delta^{-1}$.

  We get a group generated by a complex reflection
  $R_1=PJ^{-1}=\delta^2$ of angle $2\pi/p$ and a regular elliptic
  element $J$ of order 3, which satisfy the relations of
  Proposition~\ref{prop:char_s4c}.
  Discreteness then implies that the group generated by $R_1$ and $J$
  is isomorphic to $\Sc(4,\overline{\sigma}_4)$. 
\end{pf}

\begin{pf} (of Theorem~\ref{thm:iso})
  A priori the group generated by $R_1$ and $J$ as in
  Proposition~\ref{prop:key} is only a subgroup of the one generated
  by $\alpha$ and $\delta$, we show that these groups are equal, using
  the relations given in the presentation~\eqref{eq:pres_naruki}.  As
  before, we write $R_2=JR_1J^{-1}$, $R_3=J^{-1}R_1J$, $P=R_1J$. We
  then have
  $$
  R_2=\delta\bar\alpha\bar\delta\alpha\delta\alpha\bar\delta,\quad R_3=\alpha\delta\alpha\bar\delta\bar\alpha,
  $$
  and this implies $R_3^{-1}R_2R_3=\alpha$. Indeed,
  $$
  R_3^{-1}R_2R_3=\alpha\delta\bar\alpha\bar\delta\bar\alpha 
  \cdot \delta\bar\alpha\bar\delta\alpha\delta\alpha\bar\delta
  \cdot \alpha\delta\alpha\bar\delta\bar\alpha
  = \alpha\delta\bar\alpha\bar\delta\bar\alpha 
  \cdot \delta\bar\alpha\bar\delta(\alpha\delta\alpha\bar\delta)^2\bar\alpha
  = \alpha\delta\bar\alpha\bar\delta\bar\alpha 
  \cdot \delta\bar\alpha\bar\delta(\delta\alpha\bar\delta\alpha)^2\bar\alpha
  =\alpha.
  $$ 
  One can also write $\delta$ in terms of $R_1$ and $J$, namely
  $\delta=P^2R_1P^{-2}R_1P^2$ (the right hand side of this equation
  may seem complicated, but it appears in previous work of the author,
  see p.~708 in~\cite{dpp2}). Indeed,
  $P^2=\delta(\delta\alpha)^6\bar\delta=\delta\bar\alpha\bar\delta^2$,
  since $\alpha\delta$ has order 7. This gives
  $$
  P^2R_1P^{-2}R_1P^2=\delta\bar\alpha\bar\delta^2\cdot \delta^2\cdot
  \delta^2\alpha\bar\delta\cdot \delta^2\cdot
  \delta\bar\alpha\bar\delta^2
  = \delta\bar\alpha(\delta^2\alpha\delta^2)\bar\alpha\bar\delta^2
  = \delta\bar\alpha(\alpha\delta^2\alpha)\bar\alpha\bar\delta^2
  =\delta
  $$
\end{pf}

For completeness, we mention the following.
\begin{prop}\label{prop:infty}
  The lattice $G^{(\infty)}$ is an arithmetic group commensurable with
  $PU(2,1,\mathcal{O}_7)$, where $\mathcal{O}_7$ is the ring of
  integers in $Q(i\sqrt{7})$.
\end{prop}
\begin{pf}
  Proposition~\ref{prop:key} shows that $G^{(\infty)}$ is generated by
  a parabolic element $R_1$ and a regular elliptic element $J$ such
  that the relations in Proposition~\ref{prop:char_s4c} hold, i.e. it is
  isomorphic to $\Sc(\infty,\overline{\sigma}_4)$.
  Using the (obvious extension to the case $p=\infty$ of the)
  description of the groups $\Sc(p,\tau)$ given in section~2.5
  of~\cite{parkerpaupert}, we may write
  $$
  R_1=\left(\begin{matrix}
    1 & \frac{-1-i\sqrt{7}}{2} & \frac{1-i\sqrt{7}}{2}\\
    0 & 1 & 0\\
    0 & 0 & 1
    \end{matrix}\right),\quad
  J = \left(\begin{matrix}
    0 & 0 & 1\\
    1 & 0 & 0\\
    0 & 1 & 0
    \end{matrix}\right), \quad
  H = \left(\begin{matrix}
    0 & \mu & \bar\mu\\
    \bar\mu & 0 & \mu\\
    \mu & \bar\mu & 0
    \end{matrix}\right),
  $$ 
  where $\mu=-i\sqrt{7}\overline{\sigma}_4=(-7+i\sqrt{7})/2$.  This
  exhibits $\Sc(\infty,\overline{\sigma}_4)$ as a subgroup of
  $U(H,\mathcal{O}_7)$, where $H$ is a form with entries in
  $\mathcal{O}_7$. Since $\Q(i\sqrt{7})$ is a quadratic number field,
  $\mathcal{O}_7$ is discrete and $U(H,\mathcal{O}_7)$ is arithmetic.

  It follows from the classification of arithmetic subgroups of
  $PU(2,1)$ that $G^{(\infty)}$ is commensurable to
  $PU(2,1,\mathcal{O}_7)$ (see e.g.~\cite{emerystover}).
\end{pf}

\section{Congruence subgroups}\label{sec:congruence}

As mentioned in the introduction, the lattices $G^{(p)}$ for
$p=4,6,8,\infty$ appear in work of Barthel, Hirzebruch and H\"ofer,
see~\cite{bhh}. The reason why only these values of $p$ appear there
is that Barthel, Hirzebruch and H\"ofer worked on the level of $\P^2$,
not on the level of the finite quotient $\quotr{\P^2}{G_{24}}=\P(2,3,7)$. Since
the corresponding quotient map $\P^2\rightarrow \quotr{\P^2}{G_{24}}$ branches
with order 2 on the mirrors of reflections in $G_{24}$, the orbifold
structures we consider on $\quotr{\P^2}{G_{24}}$ (possibly blown-up) do not
lift to an orbifold structure on $\P^2$ (possibly blown-up). Indeed,
generic points on the Klein discriminant curve, which have integer
multiplicy $p$, would lift to points with multiplicity $p/2$ in $\P^2$
(which is in general not an integer). Similarly, for $p=12$, the curve
$E$ has odd multiplicity, see Table~\ref{tab:pairs}, and the orbifold
structure does not lift to $\P^2$ either.

The difference between even or odd weights is closely related to the
distinction between the INT and the $\Sigma$-INT condition for the
hypergeometric monodromy groups of Deligne-Mostow,
see~\cite{mostowIHES}.

In this section, we interpret the coverings corresponding to
pulling-back the orbifold structure via $\P^2\rightarrow quot{\P^2}{G_{24}}$
for $p=4,6,8,\infty$ as explicit congruence subgroups.

 In order to get explicit linear groups $\widetilde{\Gamma}_p$ in
 $U(2,1)$ (rather than the projectivized $\PU(2,1)$), we use specific
 matrices $R_1,R_2,R_3$ as generators of the sporadic group
 $\widetilde{\Gamma}_p$, namely:
\begin{equation}\label{eq:r123}
R_1=\left(\begin{matrix}
a & \tau & -\overline{\tau}\\
0 & 1 & 0\\
0 & 0 & 1
\end{matrix}\right),\quad
R_2=\left(\begin{matrix}
1 & 0 & 0\\
-a\overline{\tau} & a & \tau\\
0 & 0 & 1
\end{matrix}\right),\quad
R_3=\left(\begin{matrix}
1 & 0 & 0\\
0 & 1 & 0\\
a\tau & -a\overline{\tau} & a\\
\end{matrix}\right).
\end{equation}
where $\tau=-\frac{1+i\sqrt{7}}{2}$, $a=e^{2\pi i/p}$. Note that in
$\PU(2,1)$, we have $J=(R_1 R_2 R_3 R_1 R_2 R_3 R_1)^{-1}$, so the
matrices in equation~\eqref{eq:r123} generate the same group as $R_1$
and $J$.

We denote by $\phi_p^{(n)}:\widetilde{\mathcal{S}}_p\rightarrow
\mathbb{F}_{n^r}$ reduction modulo some prime factor of $n$ in the
field $\K=\Q(a,\tau)$. The kernel
$\widetilde{\Gamma}_p(n)=Ker(\phi_p^{(n)})$ is a congruence subgroup
of $\widetilde{\mathcal{S}}_p$. 
We
will prove the following.
\begin{thm}\label{thm:congruence}
  \begin{enumerate}
    \item $Im(\phi_4^{(2)})$, $Im(\phi_8^{(2)})$,
      $Im(\phi_\infty^{(2)})$ are all isomorphic to
      $GL_3(\mathbb{F}_2)$, which is the unique simple group of order
      168.
    \item $Im(\phi_6^{(3)})$ is a subgroup of order 336 in
      $GL_3(\F_9)$, with center of order 2. The quotient of this subgroup
      by its center is the simple group of order 168.
  \end{enumerate}
\end{thm}
In particular, the lattices $\Gamma_p$ have normal subgroups of index
168 for $p=4,6,8,\infty$. For $p=3,5,12$, there is no such normal
subgroup, as can be verified fairly easily using computational group
theory software.

\begin{pf} We give a proof for each relevant value of $p$.\\
\noindent ${\bf p=4}$

Here and in what follows, we write $\tau=-(1+i\sqrt{7})/2$.  Consider
the field $\K=\Q(i\tau)$ and the ideal generated by
$\zeta=i-\tau$. Note that $\zeta$ is a factor of 2 in
$\mathcal{O}_{\K}$, since $(i-\tau)(i-\overline{\tau})=1+i$, which
gives
$$
(i-\tau)^2(i-\overline{\tau})^2=2i.
$$ 
An integral basis for $K$ is given by $1,i(\tau-1),2+\tau,i(\tau-2)$.
One easily checks that $\mathcal{O}_{\K}/(\zeta)$ is a field with two
elements, which we simply denote by $\F_2$, and $i$ and $\tau$ reduce
to 1, while $\overline{\tau}$ reduces to $0$.

In particular, the matrices $R_1$, $R_2$ and $R_3$ from
equation~\eqref{eq:r123} reduce to
\begin{equation}\label{eq:r123-4}
\Rb_1=\left(\begin{matrix}
1 & 1 & 0\\
0 & 1 & 0\\
0 & 0 & 1
\end{matrix}\right),\quad
\Rb_2=\left(\begin{matrix}
1 & 0 & 0\\
0 & 1 & 1\\
0 & 0 & 1
\end{matrix}\right),\quad
\Rb_3=\left(\begin{matrix}
1 & 0 & 0\\
0 & 1 & 0\\
1 & 0 & 1\\
\end{matrix}\right).
\end{equation}
It is well known that $GL_3(\F_2)\cong PSL_2(\F_7)$, and that it is
the unique simple group of order 168, see the Atlas of Finite
Groups~\cite{atlas} for instance. For an elementary proof, as well as
an overview of other proofs in the literature, see
also~\cite{brownloehr}.

Using Gaussian elimination, it is easy to see that $GL_3(\F_2)$ is
generated by the above three matrices together with two permutation
matrices $T_{12}$ and $T_{23}$ corresponding to a transposition of two
standard basis vectors.

One can recover such transpositions by verifying that
$$
\Rb_1(\Rb_2\Rb_3)^2\Rb_1=T_{12},\quad \Rb_2(\Rb_3\Rb_1)^2\Rb_2=T_{23}.
$$

\noindent ${\bf p=6}$

For $p=6$, we consider $\K=\Q(\omega\tau)$, where
$\omega=(-1+i\sqrt{3})/2$, and the ideal generated by $\eta=\omega-1$,
which is a prime factor of $3$ in $\mathcal{O}_K$.

One checks that the residue field is $\F_9$, which we write as a
vector space over $\F_3$ generated by $1$ and $u$. For a suitable
choice of $u$, the reduction of the matrices are given by
\begin{equation}\label{eq:r123-6}
\Rb_1=\left(\begin{matrix}
2 & 2u & 2u+1\\
0 & 1 & 0\\
0 & 0 & 1
\end{matrix}\right),\quad
\Rb_2=\left(\begin{matrix}
1 & 0 & 0\\
u+2 & 2 & 2u\\
0 & 0 & 1
\end{matrix}\right),\quad
\Rb_3=\left(\begin{matrix}
1 & 0 & 0\\
0 & 1 & 0\\
u & u+2 & 2\\
\end{matrix}\right).
\end{equation}
The study of the group generated by these matrices is slighlty more
subtle than for the case $p=4$, since $GL_3(\F_9)$ is quite large. We
will use the group presentation for $G_{24}$ given by Shephard and
Todd, see p.~299 of~\cite{shephardtodd}. We define $A_1=\Rb_1$,
$A_2=\Rb_2$, $A_3=\Rb_2\Rb_3\Rb_2$, and verify that
$A_1^2=A_2^3=A_3^2=(A_1A_2)^4=(A_2A_3)^4=(A_3A_1)^3=(A_1A_2A_1A_3)^3=Id$. This
implies that $Im(\phi_6^{(3)})$ is a homomorphic image of $G_{24}$.

Note that $G_{24}$ has an index two subgroup which is the simple group
of order 168, and the simplicity of that subgroup implies that the
above homomorphism has trivial kernel.

\noindent ${\bf p=8}$

For $p=8$, the number field $\K=\Q(\zeta_8,\tau)$ has degree 8. The
minimal polynomial of a primitive generator, say
$\alpha=i\sqrt{7}+(1+i)/\sqrt{2}$ is given by
$x^8+28x^6+294x^4+1288x^2+2500$. The computations are clearly very
intricate, and they are best achieved with a computer. One gets a
prime factorization of the form $(2)=I_1^4I_2^4$, where $I_1$ is the
two generator ideal
{
\scriptsize
$$ 
(2,-19/15600\alpha^7 - 7/780\alpha^6 - 49/1950\alpha^5 -
59/312\alpha^4 - 1337/7800\alpha^3 - 1057/780\alpha^2 -
833/3900\alpha - 241/156).
$$ 
}
One checks that the corresponding residue field is $\F_2$, and that
the reduction mod $I_1$ gives the same matrices as in~\eqref{eq:r123-4}. Note
also that taking $I_2$ instead of $I_1$ would give an isomorphic
reduction, since $\K$ is a Galois extension of $\Q$.

\noindent ${\bf p=\infty}$

For $p=\infty$, the proof is essentially the same as for $p=4$, with
the simpler number field $\K=\Q(i\sqrt{7})$. Note that
$\mathcal{O}_{\K}$ is then simply the $\Z$-module generated by
$1,\tau$; we take the ideal generated by $\tau$, which is a prime
factor of 2, since $-\tau(\tau+1)=2$. The residue field is $\F_2$, and
one gets the same reduced matrices as in~\eqref{eq:r123-4}.

\end{pf}

\end{document}